\pdfoutput=1
\documentclass[12pt]{amsart}

\usepackage{shellesc}
\usepackage{pgfplots}
\input{pgfplots_custom.sty}
\usepgfplotslibrary{external}
	\tikzsetexternalprefix{}

\usepackage{fullpage}

\usepackage{amsmath}
\usepackage{amsfonts}
\usepackage{amssymb}
\usepackage{amsthm}
\usepackage{mathtools} 
\usepackage{latexsym}

\usepackage{enumerate}
\usepackage{cancel}
\usepackage{cases}
\usepackage{empheq}
\usepackage{multicol}

\usepackage{graphics} 

\usepackage{wrapfig}

\usepackage{color}

\usepackage[backgroundcolor=gray!30,linecolor=black]{todonotes}
   
\usepackage{mathrsfs}
\usepackage{fontenc} 
\usepackage{inputenc} 

\usepackage{verbatim}

\theoremstyle{plain}
\newtheorem{theorem}{Theorem}[section]

\newtheorem{lemma}[theorem]{Lemma}

\theoremstyle{definition}
\newtheorem{definition}[theorem]{Definition}

\theoremstyle{remark}
\newtheorem{remark}[theorem]{Remark}

\numberwithin{equation}{section} 
\numberwithin{figure}{section}   

\usepackage{url}

\usepackage{algorithm}
\usepackage{algpseudocode}

\newcommand{\vect}[1]{\mathbf{#1}}

\newcommand{\ba}{\vect{a}}
\newcommand{\bb}{\vect{b}}

\newcommand{\bD}{\vect{D}}
\newcommand{\bu}{\vect{u}}
\newcommand{\bU}{\vect{U}}
\newcommand{\bv}{\vect{v}}
\newcommand{\bV}{\vect{V}}
\newcommand{\bw}{\vect{w}}

\newcommand{\bx}{\vect{x}}

\newcommand{\bg}{\vect{g}}

\newcommand{\field}[1]{\mathbb{#1}}

\newcommand{\nR}{\field{R}}

\usepackage{mathrsfs}
\DeclareMathAlphabet{\mathpzc}{OT1}{pzc}{m}{it}

\newcommand{\RE}{\text{Re}}

\newcommand{\cnj}[1]{\overline{#1}}

\newcommand{\ip}[2]{\left<#1,#2\right>}

\newcounter{my_counter}
\title[Data Assimilation: Parameter Recovery, Sensitivity]{Parameter Recovery and Sensitivity Analysis for the 2D Navier-Stokes Equations Via Continuous Data Assimilation}

\date{\today}

\author{Elizabeth Carlson}
\address[Elizabeth Carlson]{Department of Mathematics, 
                University of Nebraska--Lincoln,
        Lincoln, NE 68588-0130, USA}
\email[Elizabeth Carlson]{elizabeth.carlson@huskers.unl.edu}
\author{Joshua Hudson}
\address[Joshua Hudson]{Johns Hopkins University Applied Physics Laboratory,
        11100 Johns Hopkins Road, 
	Laurel, MD 20723-6099, USA}
\email[Joshua Hudson]{Joshua.Hudson@jhuapl.edu}
\author{Adam Larios}
\address[Adam Larios]{Department of Mathematics, 
                University of Nebraska--Lincoln,
        Lincoln, NE 68588-0130, USA}
\email[Adam Larios]{alarios@unl.edu}

\keywords{Parameter Recovery, Sensitivity Analysis, Continuous Data Assimilation, Navier-Stokes Equations, Reynolds Number}
\thanks{MSC 2010 Classification: 
34D06, 
35A01, 
35Q30, 
35Q35, 
37C50, 
76D03 
}

\begin{document}
\begin{abstract}
We study a continuous data assimilation algorithm proposed by Azouani, Olson, and Titi (AOT) in the context of an unknown Reynolds number.  We determine the large-time error between the true solution of the 2D Navier-Stokes equations and the assimilated solution due to discrepancy between an approximate Reynolds number and the physical Reynolds number.  Additionally, we develop an algorithm that can be run in tandem with the AOT algorithm to recover both the true solution and the Reynolds number (or equivalently the true viscosity) using only spatially discrete velocity measurements.   The algorithm we propose involves changing the viscosity mid-simulation.  Therefore, we also examine the sensitivity of the equations with respect to the Reynolds number.  We prove that a sequence of difference quotients with respect to the Reynolds number converges to the unique solution of the sensitivity equations for both the 2D Navier-Stokes equations and the assimilated equations.  We also note that this appears to be the first such rigorous proof of existence and uniqueness to strong or weak solutions to the sensitivity equations for the 2D Navier-Stokes equations (in the natural case of zero initial data), and that they can be obtained as a limit of difference quotients with respect to the Reynolds number.
\end{abstract}

\maketitle
\thispagestyle{empty}


\noindent
\section{Introduction}\label{secInt}
\noindent

A major difficulty in performing accurate, practical simulations of dynamical systems is that one typically does not have complete information about the initial state of the system, nor the exact physical parameters of the system, which may be inaccurately measured, or simply unknown.  In this paper, we present an algorithm based on data assimilation that addresses both of these difficulties.  The term \textit{data assimilation} refers to a wide class of techniques for incorporating observational data into simulations to increase their accuracy.  It is especially relevant for situations in which information about the initial data is sparse.  Recently, in a paper by Azouani, Olson, and Titi \cite{Azouani_Olson_Titi_2014}, a new approach to data assimilation, which we refer to as the AOT algorithm, was proposed. This algorithm uses a feedback control term at the PDE level to penalize deviations from the observed data.  In the present work, we apply the AOT algorithm in the setting of an unknown diffusion coefficient (e.g., viscosity or $\RE^{-1}$).  Moreover, we propose a new algorithm which changes the diffusion coefficient dynamically as the simulation evolves in time, driving the parameter to its true value.

We demonstrate this parameter recovery method for estimating viscosity using the feedback control method of data assimilation proposed in \cite{Azouani_Olson_Titi_2014}, which states that given a dissipative dynamical system (of possibly infinite dimension) of the form
\[\frac{d\bu}{dt} = F(\bu)\] with missing initial data, we can instead solve the system
\begin{align*}
 \frac{d\bv}{dt} &= F(\bv) +\mu(I_h(\bu)-I_h(\bv)) \\
\bv(0) &= \bv_0,
\end{align*}
where $\mu$ is a sufficiently large positive relaxation parameter, $I_h(\bu)$ represents the observational measurements, and $\bv_0$ is arbitrarily chosen.  The function $I_h$ is a straightforward interpolant satisfying particular bounds (stated in the Preliminaries), and is often taken to be modal projection.  

Following the analysis of \cite{Azouani_Olson_Titi_2014} on the 2D incompressible Navier-Stokes equations,  analytical bounds on the large time error of $\bv$ with respect to the true solution $\bu$ are shown to be directly dependent upon the difference between our chosen Reynolds number and the true Reynolds number.  Computationally, it is observed that the term involving the error of the Reynolds numbers closely matches the error between the solutions $\bv$ and $\bu$.  Due to this fact, we develop a heuristic algorithm for computationally recovering the true viscosity and use this algorithm to simultaneously converge to the true solution $\bu$.

Since this algorithm introduces a discontinuous change in viscosity during the simulation, we want to ensure that this abrupt change does not lead to the development of shocks in the solution.  Thus, we also study the sensitivity of the systems under consideration with respect to perturbation of the diffusion parameter.  In particular, we prove that the derivative of solutions with respect to the viscosity is a well-defined object which is bounded in appropriate function spaces; additionally we prove that the corresponding sensitivity equations are globally well-posed in time in an appropriate weak sense and that weak solutions are unique.  
Sensitivity for partial differential equations has been studied formally in many contexts (see, e.g., \cite{Anderson_Newman_Whitfield_Nielsen_1999_AIAA, Breckling_Neda_Pahlevani_2018_CMA, Hyoungjin_Chongam_Rho_DongLee_1999_KSIAM, Fernandez_Moubachir_2002_MMMAS, Kouhi_Houzeaux_Cucchietti_Vazquez_2016_AIAAConf, Vemuri_Raefsky_1979_IJSS, Hamby_1994_EnvMonAssess, Davis_Pahlevani_2013, Pahlevani_2006, Borggaard_Burns_1997}).  In \cite{Stanley_Stewart_2002}, it was shown formally that the sensitivity equations for the steady-state 2D Navier-Stokes equations are globally well-posed.  Additionally, rigorous results on the existence of derivatives of solutions to generic linear and nonlinear differential equations with respect to parameters were proven in \cite{Brewer_1982_JMathAnal, Gibson_Clark_1977} via semigroup theory.  After the preparation of this manuscript, it also came to our attention that some analysis for the sensitivity equations has been carried out in the slightly more general context of a large eddy simulation (LES) model of
the 2D Navier-Stokes equations in an unpublished PhD thesis \cite{Pahlevani_2004}.  In particular, a formal proof of the global existence and
uniqueness of the equations was given, based on formal energy
estimates.  Here, for the convenience of the reader, we provide a fully
rigorous proof (although not in the more general LES setting) of the
global well-posedness of the sensitivity equations for the 2D
Navier-Stokes equations.  Moreover, we prove that difference quotients
of solutions corresponding to different viscosities converge, in an
appropriate sense, to solutions of the sensitivity equations.  In this paper, we instead rigorously prove the existence of unique weak solutions with zero initial data to the associated sensitivity equations specifically for the 2D Navier-Stokes equations using the limit of difference quotients corresponding to different viscosities.  


Our error estimates in this work are also relevant to the setting of subgrid scale data.  In real-world settings, simulations are often underresolved; in particular, it is not always possible to run simulations with the physical Reynolds number (see, e.g.,  \cite{Lesieur_Metais_Comte_2005_LES_book,Sagaut_2006_LES_book,Berselli_Iliescu_Layton_2006_book}, and the references therein).  The error estimates we prove in this paper indicate that one may simulate flows using the AOT algorithm with a Reynolds number which is, e.g., \textit{smaller} than the true Reynolds number, and be assured that deviations from the true solution are controlled (in the $L^2$ and $H^1$ norms) by the difference in the (inverse) Reynolds numbers.

We note that classical data assimilation is largely focused on statistical optimization approaches utilizing the Kalman filter \cite{Kalman_1960_JBE} or 3D/4D-Var methods, and variations of these techniques (see, e.g., \cite{Daley_1993_atmospheric_book,Kalnay_2003_DA_book,Law_Stuart_Zygalakis_2015_book,Lewis_Lakshmivarahan_2008}, and the references therein).  
The AOT algorithm (which is also called \textit{continuous data assimilation} or CDA in the literature), differs markedly from the Kalman filter approach.  Instead of employing statistical tools at the numerical level, AOT data assimilation arises at the PDE level via a feedback-control term which penalizes deviations from interpolations of observable data.  This interpolation is a key difference between the AOT method and the so-called nudging or Newtonian relaxation methods introduced in \cite{Anthes_1974_JAS,Hoke_Anthes_1976_MWR}, as it allows for significantly more sparse initial data. For an overview of nudging methods, see, e.g., \cite{Lakshmivarahan_Lewis_2013} .  We mention that a method that shares some features with the AOT algorithm was introduced in \cite{Blomker_Law_Stuart_Zygalakis_2013_NL} in the context of stochastic differential equations. 
The AOT algorithm and its extensions have been the subject of much recent theoretical work; see, e.g., \cite{Albanez_Nussenzveig_Lopes_Titi_2016,
Bessaih_Olson_Titi_2015,
Biswas_Foias_Monaini_Titi_2018downscaling,
Biswas_Martinez_2017,
Celik_Olson_Titi_2018,
Farhat_Jolly_Titi_2015,
Farhat_Lunasin_Titi_2016abridged,
Farhat_Lunasin_Titi_2016benard,
Farhat_Lunasin_Titi_2016_Charney,
Farhat_Lunasin_Titi_2017_Horizontal,
Foias_Mondaini_Titi_2016,
Foyash_Dzholli_Kravchenko_Titi_2014,
GarciaArchilla_Novo_Titi_2018,
GlattHoltz_Kukavica_Vicol_2014,
Ibdah_Mondaini_Titi_2018uniform,
Jolly_Martinez_Olson_Titi_2018_blurred_SQG,
Jolly_Martinez_Titi_2017,
Jolly_Sadigov_Titi_2015,
Larios_Pei_2018_NSV_DA,
Markowich_Titi_Trabelsi_2016_Darcy,
Mondaini_Titi_2018_SIAM_NA,
Pei_2018,
Rebholz_Zerfas_2018_alg_nudge}.  
Computational trials of the AOT algorithm and its variants were carried out on a wide variety of equations in several recent works, including \cite{Desamsetti_Dasari_Langodan_Knio_Hoteit_Titi_2018_WRF,Gesho_Olson_Titi_2015,Larios_Rebholz_Zerfas_2018,Leoni_Mazzino_Biferale_2018,Altaf_Titi_Knio_Zhao_Mc_Cabe_Hoteit_2015,Lunasin_Titi_2015,Larios_Pei_2017_KSE_DA_NL}. 
We also mention an upcoming work \cite{Farhat_GlattHoltz_Martinez_McQuarrie_Whitehead_2018}, currently a preprint, which explores some similar ideas contained in this paper in the context of continuous data assimilation for the {R}ayleigh-{B}\'enard convection equations with unknown Prandtl number, although parameter recovery is not explored in that work.

The paper is organized as follows: in Section~\ref{secPre}, we describe the mathematical framework for the problems we consider. In Section~\ref{secNeatSection} we consider the AOT data assimilation algorithm and show that with only an approximation of the true viscosity, the reference solution can still be recovered using data assimilation. We analyze the Navier-Stokes equations with periodic boundary conditions, but our techniques can be extended to other boundary conditions and other dissipative systems. In Section~\ref{secSensitivity} we prove rigorous bounds on the sensitivity of the data assimilation approximation on the approximate viscosity used. In Section~\ref{secNumerics}, we provide numerical evidence that illustrates the effectiveness of the algorithm, as well as the practical performance we might expect in a typical flow. In addition, in Section~\ref{sec:param-est}, motivated by our rigorous results we consider the inverse problem of parameter recovery, and derive an algorithm to recover the true viscosity using data assimilation. 

\section{Preliminaries}\label{secPre}
In this section, we state the theorems and other preliminaries needed to solve the incompressible Navier-Stokes equations and the associated modified equations utilizing the AOT algorithm \cite{Azouani_Olson_Titi_2014}.  The statements given in the section without proof  are standard, and proofs can be found, e.g., in, e.g.,  \cite{Constantin_Foias_1988,Foias_Manley_Rosa_Temam_2001,Robinson_2001,Temam_2001_Th_Num,Temam_1995_Fun_Anal}. We consider the incompressible Navier-Stokes equations in dimensionless form on a spatial domain $\Omega$,

\begin{subequations} \label{NSEpre}
\begin{alignat}{2}
\label{NSE_mo_pre}
\partial_t\bu + (\bu\cdot\nabla)\bu &=-\nabla p + \RE_1^{-1} \triangle\bu + \mathbf{f},
\qquad&& \text{in }\Omega\times[0,T],\\
\label{NSE_div_pre}
\nabla \cdot \bu &=0,
\qquad&& \text{in }\Omega\times[0,T],\\
\label{NSE_IC_pre}
\bu(\bx,0)&=\bu_0(\bx),
\qquad&& \text{in }\Omega.
\end{alignat}
\end{subequations}
where $\RE_1 = \frac{UL}{\nu_1}$ is the dimensionless Reynolds number based on the kinematic viscosity $\nu_1 > 0$, a typical length scale $L$, and typical velocity $U$.

We take the spatial domain, $\Omega$, to be the torus, i.e. $\Omega = \mathbb{T}^2 = \mathbb{R}^2/\mathbb{Z}^2$, which is an open, bounded, and connected domain with $C^2$ boundary.  As is customary, we define the space \[\mathcal{V}:= \{f: \Omega \to \mathbb{R}^2 \; | \;f \in \dot{C}_p^\infty(\mathbb{T}^2)\},\] and subsequently the spaces $H := \cnj{\mathcal{V}}$ in $L^2(\Omega; \nR^2)$ and $V := \cnj{\mathcal{V}}$ in $H^1(\Omega; \nR^2)$.  $H$ and $V$ are subspaces of $L^2(\Omega; \nR^2)$ and $H^1(\Omega; \nR^2)$, respectively, and hence are Hilbert spaces with the inner products defined as
\begin{alignat}{2}
(\bu,\bv) = \int_{\mathbb{T}^2} \bu \cdot \bv  \; d\bx \qquad &&  ((\bu,\bv))= \sum\limits_{i,j =1}^2 \int_{\mathbb{T}^2} \frac{\partial u_i}{\partial x_j}\frac{\partial v_i}{\partial x_j} \;d\bx, \notag
\end{alignat}
with corresponding norms $|\bu| = \sqrt{(\bu,\bu)}$ and $\|\bu\| = \sqrt{((\bu,\bu))}$.  Due to the mean-zero condition, the following Poincar\'e inequalities hold.
$$\lambda_1\|\bu\|_{L^2}^2\leq\|\nabla\bu\|_{L^2}^2 \text{\quad for\quad} \bu\in V,$$
$$\lambda_1\|\nabla\bu\|_{L^2}^2\leq\|A\bu\|_{L^2}^2 \text{\quad for\quad} \bu\in D(A).$$
Thus, $|\nabla\bu|$ and $\|\bu\|$ are equivalent norms on $V$.  In 2D, the following Brezis-Gallouet inequality, proven in \cite{Brezis_Gallouet_1980}, also holds for all $\bu\in\mathcal{D}(A)$
\begin{align}
\|\bu\|_{L^\infty} \leq c\|\bu\| \left\{1+\log \frac{|A\bu|^2}{4\pi^2\|\bu\|^2}\right\}.
\end{align}

We can equivalently consider the Leray projection of the \eqref{NSEpre}, where $P_\sigma$ is the orthogonal projection from $L^2(\Omega)$ onto $H$.  As in \cite{Azouani_Olson_Titi_2014}, we define the Stokes operator $A$ and the bilinear term $B: V\times V \to V^*$ as the continuous extensions of the operators $A$ and $B$ defined on $\mathcal{V} \times \mathcal{V}$ as

\begin{alignat}{3}
A\bu = -P_\sigma \triangle \bu & \qquad \text{ and } \qquad & B(u,v) = P_\sigma(\bu \cdot \nabla \bv),\notag
\end{alignat}
and we define the domain of $A$ to be $\mathcal{D}(A) := \{u\in V: Au\in H\}$.  Also note that $A$ is a linear self-adjoint and positive definite operator with a compact inverse, so there exists a complete orthonormal set of eigenfunctions $w_i$ in $H$ such that $Aw_i = \lambda_iw_i$, with the eigenvalues strictly positive and monotonically increasing.  

We note that the bilinear operator, $B$, has the property
\begin{align}\label{bilinear_symmetry}
\ip{B(\bu,\bv)}{\bw} = -\ip{B(\bu,\bw)}{\bv},
\end{align}
for all $\bu,\bv,\bw \in V$.  This implies $B$ also satisfies
\begin{align}\label{Borth}
\ip{B(\bu,\bw)}{\bw} = 0,
\end{align}
for all $\bu,\bv,\bw \in V$.  Moreover,  the following inequalities hold:
\begin{align}
\ip{B(\bu,\bv)}{\bw}|&\leq \|\bu\|_{L^\infty(\Omega)}\|\bv\||\bw|  &\text{ for } \bu \in L^\infty(\Omega), \bv \in V, \bw \in H
\label{BINsimple} \\
|\ip{B(\bu,\bv)}{\bw}| &\leq c |\bu|^{1/2}\|\bu\|^{1/2} \|\bv\||\bw|^{1/2}\|\bw\|^{1/2} &\text{ for } \bu,\bv,\bw \in V, \label{BIN} \\
|(B(\bu,\bv),\bw)| &\leq c|\bu|^{1/2}\|\bu\|^{1/2}\|\bv\|^{1/2}|A\bv|^{1/2}|\bw|  &\text{ for } \bu\in V, \bv \in \mathcal{D}(A), \bw \in H \label{BIN2} \\
|(B(\bu,\bv),\bw)| &\leq c |\bu|^{1/2}|A\bu|^{1/2} \|\bv\||\bw| &\text{ for } \bu \in \mathcal{D}(A), \bv \in V, \bw \in H.\end{align}

Due to the periodic boundary conditions, it also holds (in 2D) that
\begin{align}\label{bilinear identity}
(B(\bw,\bw),A\bw)=0 \quad \text{ for every } \quad \bw \in \mathcal{D}(A).
\end{align}
Therefore, for $\bu, \bw \in \mathcal{D}(A)$,
\begin{align}\label{lastbilinear}
(B(\bu,\bw),A\bw)+(B(\bw,\bu),A\bw) = - (B(\bw,\bw),A\bu).
\end{align}

Additionally, we have the following properties of the bilinear term Lemma \ref{bilinear_unif_bd} and \ref{bilinear_wk_conv}, which we prove using similar strategies as in \cite{Robinson_2001}.

\begin{lemma}\label{bilinear_unif_bd}
 Suppose $\{\ba_n\}_{n\in\mathbb{N}}$ and $\{\bb_n\}_{n\in\mathbb{N}}$ are uniformly bounded sequences in $L^2(0,T;V) \cap L^\infty(0,T;H)$. Then $\|B(\ba_n,\bb_n)\|_{L^2(0,T;V*)}$ is uniformly bounded in $n$.  Moreover, if $\{\ba_n\}_{n\in\mathbb{N}}$ and $\{\bb_n\}_{n\in\mathbb{N}}$ are uniformly bounded in $L^2(0,T;\mathcal{D}(A)) \cap L^\infty(0,T;V)$, then $\|B(\ba_n,\bb_n)\|_{L^2(0,T;H)}$ is uniformly bounded in $n$.
\end{lemma}
\begin{proof}
By the definition of the dual norm and \eqref{bilinear_symmetry}, 
 \begin{align*}
  \|B(\ba_n,\bb_n)\|_{V^*} &= \sup_{\stackrel{\bw \in V}{\|\bw\|=1}} |(B(\ba_n,\bb_n),\bw)|\\
  &= \sup_{\stackrel{\bw \in V}{\|\bw\|=1}} |(B(\ba_n,\bw),\bb_n)|,
 \end{align*}
and applying \eqref{BIN} we obtain
\begin{align*}
  \|B(\ba_n,\bb_n)\|_{V^*} &= \sup_{\stackrel{\bw \in V}{\|\bw\|=1}} |(B(\ba_n,\bw),\bb_n)|\\
  &\leq \sup_{\stackrel{\bw \in V}{\|\bw\|=1}} k |\ba_n|^{1/2}\|\ba_n\|^{1/2}|\bb_n|^{1/2}\|\bb_n\|^{1/2}\|\bw\| \\
  &= k |\ba_n|^{1/2}\|\ba_n\|^{1/2}|\bb_n|^{1/2}\|\bb_n\|^{1/2}.
 \end{align*}
  Using H{\"o}lder's inequality,
  \begin{align*}
   \|B(\ba_n,\bb_n)\|_{L^2(0,T;V^*)} &\leq \int_0^T \|B(\ba_n(s),\bb_n(s)\|_{V^*}^2 ds \\
   &\leq \int_0^T k |\ba_n|\|\ba_n\||\bb_n|\|\bb_n\| ds \\
   &\leq k \|\ba_n\|_{L^\infty(0,T;H)}\|\bb_n\|_{L^\infty(0,T;H)} \int_0^T \|\ba_n(s)\|\|\bb_n(s)\| ds \\
   &\leq k \|\ba_n\|_{L^\infty(0,T;H)}\|\bb_n\|_{L^\infty(0,T;H)}\|\ba_n\|_{L^2(0,T;V)}\|\bb_n\|_{L^2(0,T;V)}.
  \end{align*}
  Hence, since $\{\ba_n\}_{n\in\mathbb{N}}$ and $\{\bb_n\}_{n\in\mathbb{N}}$ are uniformly bounded  in $L^2(0,T;V) \cap L^\infty(0,T;H)$, it follows that 
$\|B(\ba_n,\bb_n)\|_{L^2(0,T;V^*)}$ is uniformly bounded in $n$.

Next, suppose $\{\ba_n\}_{n\in\mathbb{N}}$ and $\{\bb_n\}_{n\in\mathbb{N}}$ are uniformly bounded sequences in $L^2(0,T;\mathcal{D}(A)) \cap L^\infty(0,T;V)$.  Then by definition,
\begin{align*}
\|B(\ba_n,\bb_n)\|_{L^2(0,T;H)} &= \int_0^T |B(\ba_n,\bb_n)|^2 dt \\
&\leq \int_0^T |\ba_n|^2\|\bb_n\|^2 dt \\
&\leq \|\bb_n\|^2_{L^\infty(0,T;V)} \int_0^T |\ba_n|^2 dt\\
&\leq \frac{1}{\lambda_1^2}\|\bb_n\|^2_{L^\infty(0,T;V)} \int_0^T |A\ba_n|^2 dt\\
&= \frac{1}{\lambda_1^2}\|\bb_n\|^2_{L^\infty(0,T;V)} \|\ba_n\|_{L^2(0,T;\mathcal{D}(A))},
\end{align*}
which implies $\|B(\ba_n,\bb_n)\|_{L^2(0,T;H)}$ is uniformly bounded in $n$.
\end{proof}

\begin{lemma}\label{bilinear_wk_conv_v*}
  Let $\ba, \bb \in L^2(0,T;H)$, and suppose that $\ba_n \to \ba$ and $\bb_n \to \bb$ strongly in $L^2(0,T;H)$.  Suppose also that the sequences $\{\ba_n\}$ and $\{\bb_n\}$ are bounded above uniformly in $n$ in $L^\infty(0,T;H)$.  Then $B(\ba_n,\bb_n) \stackrel{*}{\rightharpoonup} B(\ba,\bb)$ in $L^2(0,T;V^*)$.
\end{lemma}
\begin{proof}
Take $\bw \in C^1(0,T;C^1(\Omega))$; then, with
\begin{align*}
 \int_0^T (B(\ba_n,\bb_n),\bw) dt &= -\int_0^T (B(\ba_n,\bw),\bb_n) dt \\
 &= -\sum\limits_{i,j=1}^2 \int_0^T\int_\Omega (a_n)_i(D_iw_j)(b_n)_j dx dt.
\end{align*}
This implies that
\begin{align*}
 &\quad\int_0^T (B(\ba_n,\bb_n),\bw) - (B(\ba,\bb),\bw) dt \\
 &= \int_0^T -(B(\ba_n,\bw),\bb_n) + (B(\ba,\bw),\bb) dt \\
 &= \sum\limits_{i,j=1}^2 \int_0^T\int_\Omega -(a_n)_i(D_iw_j)(b_n)_j+(a)_i(D_iw_j)(b)_j dx dt\\
 &= \sum\limits_{i,j=1}^2 \int_0^T\int_\Omega ((a)_i - (a_n)_i)(D_iw_j)(b_n)_j \\
 & \phantom{========} + ((b)_j - (b_n)_j)(D_iw_j)(a)_i dx dt. \notag
\end{align*}
Since $\ba_n \to \ba$ in $L^2(0,T;H)$, $\bb_n \to \bb$ in $L^2(0,T;H)$, and the sequences are bounded above uniformly in $L^\infty(0,T;H)$, we follow the argument in \cite{Robinson_2001} to obtain  \[\int_0^T (B(\ba_n,\bb_n),\bw) - (B(\ba,\bb),\bw) dt \to 0\]
as $n \to \infty$, and therefore by the density of $C^1(0,T;C^1(\Omega))$ in $L^2(0,T;V)$, $B(\ba_n,\bb_n) \stackrel{*}{\rightharpoonup} B(\ba,\bb)$.
\end{proof}

\begin{lemma}\label{bilinear_wk_conv}
If $\ba, \bb \in L^2(0,T;V)$, $\ba_n \to \ba$ and $\bb_n \to \bb$ strongly in $L^2(0,T;V)$, and $\{\ba_n\}$ and $\{\bb_n\}$ are bounded above uniformly in $n$ in $L^2(0,T;\mathcal{D}(A))$, then $B(\ba_n,\bb_n) \rightharpoonup B(\ba,\bb)$ in $L^2(0,T;H)$.
\end{lemma}
\begin{proof}
Take $\bw \in C(0,T;H)$; then
\begin{align*}
 &\quad\int_0^T (B(\ba_n,\bb_n),\bw) - (B(\ba,\bb),\bw) dt \\
 &= \int_0^T (B(\ba_n-\ba,\bb_n),\bw) + (B(\ba,\bb_n - \bb),\bw) dt \\
 &\leq \int_0^T |(B(\ba_n-\ba,\bb_n),\bw)| dt + \int_0^T |(B(\ba,\bb_n - \bb),\bw)| dt
\end{align*}
Applying \eqref{BINsimple}, \eqref{BIN2}, and Poincar{\'e}'s inequality we obtain
\begin{align*}
&\quad\int_0^T (B(\ba_n,\bb_n),\bw) - (B(\ba,\bb),\bw) dt \\
 &\leq c\lambda_1^{-1/2} \int_0^T \|\ba_n-\ba\| |A\bD^n| |w| dt + c\int_0^T \|\ba\|_{L^\infty(\Omega)}\|\bb_n-\bb\||\bw| dt.
\end{align*}
Applying Agmon's inequality, 
\begin{align*}
&\quad 
\int_0^T (B(\ba_n,\bb_n),\bw) - (B(\ba,\bb),\bw) dt \\
 &\leq c\lambda_1^{-1/2} \int_0^T \|\ba_n-\ba\| |A\bD^n| |w| dt + c\int_0^T |\ba|^{1/2}|A\ba|^{1/2}\|\bb_n-\bb\||\bw| dt \\
 &\leq c\lambda_1^{-1/2} \int_0^T \|\ba_n-\ba\| |A\bD^n| |w| dt + c\lambda_1^{1/2}\int_0^T |A\ba|\|\bb_n-\bb\||\bw| dt\\
 &\leq c\lambda_1^{-1/2} \|\ba_n-\ba\|_{L^2(0,T;V)}\|\bw\|_{L^\infty(0,T;H)} \|\bb_n\|_{L^2(0,T;\mathcal{D}(A))} \\
 &\phantom{=} + c\lambda_1^{-1/2} \|\ba\|_{L^2(0,T;\mathcal{D}(A))} \|\bw\|_{L^\infty(0,T;H)}\|\bb_n-\bb\|_{L^2(0,T;V)}
\end{align*}
Since $\ba_n \to \ba$ in $L^2(0,T;V)$, $\bb_n \to \bb$ in $L^2(0,T;V)$, the sequences are bounded above uniformly in $L^2(0,T;\mathcal{D}(A))$, and $\bw$ is continuous in time, then \[\int_0^T (B(\ba_n,\bb_n),\bw) - (B(\ba,\bb),\bw) dt \to 0\]
as $n \to \infty$, and therefore by the density of $C(0,T;H)$ in $L^2(0,T;H)$, $B(\ba_n,\bb_n) \rightharpoonup B(\ba,\bb)$ in $L^2(0,T;H)$.
\end{proof}

Without loss of generality, we will assume $\mathbf{f} \in L^\infty(0,T;H)$ so that $P_\sigma f = f$.
Thus, we may rewrite \eqref{NSEpre} as
\begin{subequations} \label{NSE}
\begin{alignat}{2}
\label{NSE_mo} 
\frac{d}{dt}\bu +B(\bu,\bu) &=  \RE_1^{-1} A\bu + \mathbf{f},
\qquad&& \text{in }\Omega\times[0,T],\\
\label{NSE_IC}
\bu(\bx,0)&=\bu_0(\bx),
\qquad&& \text{in }\Omega.
\end{alignat}
\end{subequations}
The pressure term can be recovered using de Rham's theorem \cite{Temam_2001_Th_Num, Foias_Manley_Rosa_Temam_2001}, a corollary of which is that
\begin{align}
\bg = \nabla p \text{ with $p$ a distribution if and only if } \ip{\bg}{\text{\textbf{h}}} = 0 \text{ for all \textbf{h}$ \in \mathcal{V}$. }
\end{align}

For a given force $f$ and some initial data $\bu_0$, it is classical that a unique global solution $\bu$ of $\eqref{NSE}$ will exist. However, we don't expect to know $\bu_0$ exactly, and so cannot compute $\bu(t)$ from $\eqref{NSE}$; rather, we consider the case that measurement data is collected on $\bu(t)$ over the time interval $[0,T]$, sufficient for the interpolation operator $I_h$ to construct the interpolation $I_h(\bu(t))$ on $[0,T]$. From here, we can define a new system, dubbed the data assimilation system, by introducing a feedback control (nudging term) via $I_h$ into \eqref{NSE} (or \eqref{NSEpre}), as is done in \cite{Azouani_Olson_Titi_2014}. 

We will construct our data assimilation system under the more general case of having only an approximate Reynolds number, $\RE_2$:
\begin{subequations}\label{NSEmodified}
\begin{alignat}{2}
\label{NSEmodified_mo}
\frac{d}{dt}\bv + B(\bv,\bv) &= \RE_2^{-1}A\bv + \mathbf{f} + \mu P_\sigma(I_h(\bu)-I_h(\bv)) \\
\label{NSEmodified_IC}
\bv(\bx,0)&=\bv_0(\bx).
\end{alignat}
\end{subequations}
Here, $\mu >0$ is a relaxation parameter, $\RE_2 = \frac{UL}{\nu_2}$ with $\nu_2$ a kinematic viscosity approximating $\nu_1$, and $I_h$ is a linear interpolant satisfying 
\begin{align}\label{IntIN}
\|\varphi-I_h(\varphi)\|^2_{L^2(\Omega)} \leq c_0 h^2\|\varphi\|^2_{H^1(\Omega)}
\end{align}

From \cite{Azouani_Olson_Titi_2014}, \eqref{NSEmodified} has a unique solution given either no-slip Dirichlet or periodic boundary conditions as stated in the following theorem.
\begin{theorem}
Suppose $I_h$ satisfies \eqref{IntIN} and $\mu c_0h^2 \leq \RE_2^{-1}$, where $c_0$ is the constant from \eqref{IntIN}.  Then the continuous data assimilation equations \eqref{NSEmodified} possess unique strong solutions that satisfy
\begin{alignat}{2}
\bv \in C([0,T];V) \cap L^2((0,T); D(A)) && \text{and } \; \frac{d\bv}{dt} \in L^2((0,T);H),
\end{alignat}
for any $T > 0$.  Furthermore, this solution is in $C([0,T],V)$ and depends continuously on the initial data $\bv_0$ in the $V$ norm.
\end{theorem}

For equations \eqref{NSE} and \eqref{NSEmodified}, we denote the dimensionless Grashof numbers as
\begin{align}
G_1 &= \frac{\RE_1^2}{4\pi^2} \limsup\limits_{t \to \infty} \|\mathbf{f}(t)\|_{L^2(\Omega)} \label {G1} \\ 
G_2 &= \frac{\RE_2^2}{4\pi^2} \limsup\limits_{t \to \infty} \|\mathbf{f}(t)\|_{L^2(\Omega)} \label{G2},
\end{align}
where $\lambda_1 = 4\pi^2 >0$ is the first eigenvalue of the Stokes operator.  In 2D, it is classical that \eqref{NSEpre} possesses a unique global strong solution. Furthermore, we have explicit upper bounds on the norms of the solution in $H$ and $V$ in terms of $G_1$.
\begin{theorem}\label{ubounds}
Fix $T >0$. Suppose that $\bu$ is a solution of $\eqref{NSE}$, corresponding to the initial value $\bu_0 \in V$. Then there exists a time $t_0$ which depends on $\bu_0$ such that for all $t \geq t_0$, it holds that
\begin{align}
|\bu(t)|^2 \leq 2\frac{G_1^2}{\RE_1^2} \; \text{ and } \int\limits_t^{t+T} \|\bu(\tau)\|^2 d\tau \leq 2\left(1+T\frac{4\pi^2}{\RE_1}\right) \frac{G_1^2}{\RE_1}.
\end{align}
In the case of periodic boundary conditions it also holds for all $t \geq t_0$ that
\begin{align}
\| \bu(t)\|^2 \leq 2\frac{4\pi^2G_1^2}{\RE_1^2}, \qquad \int\limits_t^{t+T} |A\bu(\tau)|^2 d\tau \leq 2\left(1+T\frac{4\pi^2}{\RE_1}\right) \frac{4\pi^2G_1^2}{\RE_1}.
\end{align}
furthermore, if $\mathbf{f} \in H$ is time-independent then
\begin{align}\label{Au L2 bound}
|A\bu(t)|^2 \leq \frac{16\pi^4c}{\RE_1^2} (1+G_1)^4.
\end{align}
\end{theorem}

To prove our main theoretical results, we will need the following corollary of the statement of the uniform Gr{\"o}nwall lemma proved in \cite{Jones_Titi_1992}.

\begin{lemma}[Generalized Uniform Gr{\"o}nwall Inequality] \label{GronwallLemma} Let $\alpha$ be a locally integrable real-valued function defined on $(0,\infty)$, satisfying the following conditions for some $0 < T< \infty$:
\begin{align*}
\liminf\limits_{t\to\infty} \int\limits_t^{t+T} \alpha(\tau) \; d\tau = \gamma > 0,
\end{align*} 

\begin{align*}
\limsup\limits_{t\to\infty}\int\limits_t^{t+T} \alpha^-(\tau) \; d\tau = \Gamma < \infty,
\end{align*}
where $\alpha^- = \max\{-\alpha,0\}$.  Furthermore, let $\beta$ be a real-valued locally integrable function defined on $(0,\infty)$, and let $\beta^+= \max\{\beta,0\}$.
Suppose that $\xi$ is an absolutely continuous non-negative function on $(0,\infty)$ such that 
\begin{align}\label{gronwallineq}
\frac{d}{dt} \xi + \alpha \xi \leq \beta \qquad \text{ a.e. on $(0,\infty)$}.
\end{align}
Then 
\begin{align*}
\xi(t) \leq \xi(t_0)\Gamma'e^{-\frac{\gamma}{2T}(t-t_0)} + \left(\sup\limits_{t\geq t_0}\int\limits_{t}^{t+T} \beta^+(\tau)\; d\tau \right) \Gamma' \frac{e^{\gamma/2}}{e-1},
\end{align*}
where $\Gamma' = e^{\Gamma +1+\gamma/2}$ and $t_0$ is chosen sufficiently large so that, for all $s \geq t_0$,
\begin{align}\label{alphaGamma}
\int\limits_s^{s+T} \alpha^-(\sigma) \; d\sigma \leq \Gamma +1
\end{align}
and
\begin{align}\label{alphagamma}
\int\limits_s^{s+T} \alpha(\sigma) \; d\sigma \geq \gamma/2.
\end{align}
\end{lemma}
We will also make use the following lemma proved in \cite{Azouani_Olson_Titi_2014}.
\begin{lemma}\label{log inequality}
Let $\phi(r) = r-\beta(1+\log r)$ where $\beta > 0$.  Then
\begin{align*}
\min\{\phi(r) : r \geq 1\} \geq -\beta\log \beta.
\end{align*}
\end{lemma}

 \section{Error of Continuous Data Assimilation to Viscosity}\label{secNeatSection}

We now present our first result. In \cite{Azouani_Olson_Titi_2014}, it was shown for the case $\RE_1 = \RE_2$, that given a strong solution $\bu$ of \eqref{NSEpre} and an interpolant $I_h$ satisfying \eqref{IntIN}, for sufficiently large $\mu$ and sufficiently small $h$, the corresponding solution $\bv$ of \eqref{NSEmodified} will converge in the $L^2$ sense to $\bu$ exponentially fast in time for any $\bv_0\in V$, (and convergence in the $H^1$ sense under stronger  smoothness assumptions). We extend this result to include the case $\RE_1 \neq \RE_2$.  In particular, we show that the $L^2$ error decays exponentially in time, down to a level which is controlled by the difference in the (inverse) Reynolds numbers.  Moreover, this level goes to zero as $\RE_2\rightarrow\RE_1$.  This means that the AOT algorithm for 2D Navier-Stokes can recover the solution approximately even when the true Reynolds number (equivalently, the true viscosity) is unknown, and that the accuracy improves as the approximation of the Reynolds number improves, and with the same order.

\begin{theorem}\label{thmMainResult1}
Let $\bu$ and $\bv$ be solutions to the systems \eqref{NSE} and \eqref{NSEmodified}, respectively, with initial data $\bu_0$, $\bv_0\in H$.  Suppose $\RE_1$, $\RE_2>0$.  Let $\mu \geq 20\pi^2c^2  \frac{\RE_2}{\RE_1^2}G_1^2$ and $h \leq \Big(\frac{1}{32\pi^2 c^2c_0 G_1^2}\frac{\RE_1^2}{\RE_2^2}\Big)^{1/2}$.  Then for any $T$ such that $\frac{\RE_1}{4\pi^2} < T< \infty$, and for a.e. $t>T$, it holds that
\begin{align*}
|\bv(t)-\bu(t)|^2 &\leq |\bv(t_0)-\bu(t_0)|^2 e^{1+\gamma/2} e^{-\frac{\gamma}{2T}(t-t_0)} + C \cdot \RE_2(\RE_2^{-1}-\RE_1^{-1})^2,
\end{align*}
where  \[C:= \frac{e^{1+\gamma}}{e-1}\left(2(1+4\pi^2T\RE_1^{-1})\frac{1}{\RE_1}G_1^2\right)\] and 
\begin{align*}
\gamma :=\liminf\limits_{t \to \infty} \int_t^{t+T} \mu - 2c^2 \RE_1\|\bu(s)\|^2 \; ds >0.
\end{align*}
In particular,
\begin{align*}
\limsup\limits_{t\to \infty} |\bv(t)-\bu(t)|
&\leq 
C \sqrt{\RE_2}|\RE_2^{-1}-\RE_1^{-1}|
=
C \frac{|\RE_2-\RE_1|}{\RE_1\sqrt{\RE_2}}.
\end{align*}
\end{theorem}
The idea of the proof is similar to the proof of the corresponding result in \cite{Azouani_Olson_Titi_2014}, except that we have an additional term to handle since we allow for the case $\RE_1 \neq\RE_2$.
\begin{proof}
We subtract \eqref{NSE_mo} from \eqref{NSEmodified_mo} to obtain
\begin{align*}
\bw_t + B(\bw,\bu)+B(\bv, \bw) = -\RE_1^{-1}A\bu + \RE_2^{-1} A\bv - \mu P_\sigma( I_h(\bw)),
\end{align*}
which can be simplified to

\begin{align}
\label{NSE_modiff}
\bw_t + B(\bw,\bu)+ B(\bv, \bw) &= (\RE_2^{-1}-\RE_1^{-1}) A \bu + \RE^{-1}_2 A \bw - \mu P_\sigma( I_h(\bw)),
\end{align}
 with initial data given by
 \begin{align} \notag
 \bw(\bx,0)&=\bw_0(\bx) := \bu_0(x)-\bv_0(x).
 \end{align}

We take the action of \eqref{NSE_modiff} on $\bw$, and utilize the Cauchy-Schwarz and Young's inequalities to obtain

\begin{align*}
&\quad
\frac{1}{2}\frac{d}{dt} |\bw|^2 + \RE_2^{-1} \|\bw\|^2
\\&= \notag
(\RE_2^{-1}-\RE_1^{-1}) (\nabla\bu,\nabla \bw) - \ip{B(\bw,\bu)}{\bw} - \mu(P_\sigma( I_h(\bw),\bw)
\\&\leq \notag
|\RE_2^{-1}-\RE_1^{-1}|\|\bu\|\|\bw\|  \\
&\phantom{=} \notag- \ip{B(\bw,\bu)}{\bw}-  \mu(P_\sigma( I_h(\bw)),\bw)
\\&\leq \notag
\frac{\RE_2|\RE_2^{-1}-\RE_1^{-1}|^2}{2}\|\bu\|^2+\frac{1}{2\RE_2}\|\bw\|^2 \\
&\phantom{=} \notag - \ip{B(\bw,\bu)}{\bw}-  \mu(P_\sigma( I_h(\bw),\bw).
\end{align*}
Using \eqref{BIN}, \eqref{IntIN}, and Young's inequality, we obtain
\begin{align*}
&\quad
\frac{1}{2}\frac{d}{dt} |\bw|^2 + \frac{1}{2\RE_2}\|\bw\|^2 
\\&\leq \notag
-\ip{B(\bw,\bu)}{\bw} + \frac{\RE_2|\RE_2^{-1}-\RE_1^{-1}|^2}{2}\|\bu\|^2 - \mu(P_\sigma( I_h(\bw)),\bw)
\\&\leq \notag
c |\bw| \|\bw\|\|\bu\| + \frac{\RE_2(\RE_2^{-1}-\RE_1^{-1})^2}{2}\|\bu\|^2
\\&\quad \notag
- \mu(P_\sigma( I_h(\bw)-\bw),\bw)-\mu\|\bw\|^2
\\&\leq \notag
c |\bw| \|\bw\|\|\bu\| +  \frac{\RE_2(\RE_2^{-1}-\RE_1^{-1})^2}{2}\|\bu\|^2
\\&\quad \notag
+ \mu \sqrt{c_0 h^2} \|\bw\||\bw| - \mu\|\bw\|^2
\\&\leq \notag
c^2\RE_2 |\bw|^2 \|\bu\|^2 +\frac{1}{4\RE_2} \|\bw\|^2
\\&\quad \notag
+   \frac{\RE_2(\RE_2^{-1}-\RE_1^{-1})^2}{2}\|\bu\|^2+ \mu \sqrt{c_0 h^2} \|\bw\||\bw| - \mu\|\bw\|^2
\\&\leq \notag
c^2\RE_2 |\bw|^2\|\bu\|^2 +\frac{1}{4\RE_2} \|\bw\|^2
\\&\quad \notag
+ \frac{\RE_2(\RE_2^{-1}-\RE_1^{-1})^2}{2}\|\bu\|^2 +\frac{\mu c_0h^2}{2}\|\bw\|^2- \frac{\mu}{2} |\bw|^2.
\end{align*}

This implies 
\begin{align}\label{ineq:before-Gronwall-0}
&\quad
\frac{1}{2} \frac{d}{dt} |\bw|^2 + \Big(\frac{1}{4\RE_2} - \frac{\mu c_0 h^2}{2}\Big)\|\bw\|^2 + \Big(\frac{\mu}{2}-c^2\RE_2\|\bu\|^2\Big)|\bw|^2
\\&\leq \notag
\frac{\RE_2(\RE_2^{-1}-\RE_1^{-1})^2}{2}\|\bu\|^2.
\end{align}

Since, by assumption, \[\mu \geq 20\pi^2c^2 \frac{\RE_2}{\RE_1^2}G_1^2\] and \[h \leq \Big(\frac{1}{8 c^2c_04\pi^2 G_1^2}\frac{\RE_1^2}{\RE_2^2}\Big)^{1/2},\]
it follows that
\begin{align}\label{ineq:before-Gronwall-1}
&\quad
\frac{1}{2} \frac{d}{dt} |\bw|^2 + \left(\frac{\mu}{2}-c^2\RE_2\|\bu\|^2\right)|\bw|^2
\leq 
\frac{\RE_2(\RE_2^{-1}-\RE_1^{-1})^2}{2}\|\bu\|^2.
\end{align}
Hence, we have an inequality of the form \eqref{gronwallineq}.  

Fix $T>0$ such that $\frac{\RE_1}{4\pi^2} < T < \infty$.  Then
\begin{align*}
\liminf\limits_{t \to \infty} \int_t^{t+T} \mu - 2c^2 \RE_1\|\bu(s)\|^2 \; ds
\\ \geq T\mu - 2c^2\RE_2(2(1+4\pi^2T\RE_1^{-1})\RE_1^{-1}G_1^2)) > 0,
\end{align*}
thanks to the assumption $\mu \geq 20\pi^2c^2\frac{\RE_2}{\RE_1^2} G_1^2$.  Define 
\begin{align*}
\gamma :=\liminf\limits_{t \to \infty} \int_t^{t+T} \mu - 2c^2 \RE_1\|\bu(s)\|^2 \; ds >0.
\end{align*}


Choose $t_0$ sufficiently large so that Theorem \ref{ubounds} holds and the inequalities \eqref{alphaGamma} and \eqref{alphagamma} hold.
Then
\begin{align*}
\Gamma:= \limsup\limits_{t\to\infty} \int\limits_{t}^{t+T} \alpha^-(\tau) \; d\tau = 0 < \infty,
\end{align*}
and we can apply Lemma \ref{GronwallLemma} to conclude that, for a.e. $t>t_0$, 
\begin{align*}
|\bw(t)|^2 &\leq |\bw(t_0)|^2 e^{1+\gamma/2}e^{-\frac{\gamma}{2T}(t-t_0)} + \left(\sup\limits_{t\geq t_0}\int\limits_{t}^{t+T} \RE_2 (\RE_2^{-1}-\RE_1^{-1})^2 \|\bu(\tau)\|^2\; d\tau \right) \frac{e^{1+\gamma}}{e-1} \\
&\leq |\bw(t_0)|^2 e^{1+\gamma/2}e^{-\frac{\gamma}{2T}(t-t_0)} + C \cdot \RE_2(\RE_2^{-1}-\RE_1^{-1})^2, \notag
\end{align*}
where $C:= \frac{e^{1+\gamma}}{e-1}\left(2(1+4\pi^2T\RE_1^{-1})\frac{1}{\RE_1}G_1^2\right)$.
%
%
%
Taking the limit supremum as $t\rightarrow0$ establishes the result.
\end{proof}

We now prove a similar result for the $H^1$ norm of the difference of the solutions, the proof of which closely follows that of \cite{Azouani_Olson_Titi_2014}, although again with an additional term to allow for $\RE_1 \neq\RE_2$.

\begin{theorem}\label{thmMainResult2}
Given the systems \eqref{NSE} and \eqref{NSEmodified} with periodic boundary conditions, and given $\mu \geq 12\pi^2\RE_1^{-1}JG_1$, with \[J :=\left[ 2c \log \left(\frac{2c^{3/2}\RE_2}{\RE_1}\right) + 4c\log(1+G_1)\right],\] and $\mu c_0 h^2 \leq \RE_2^{-1}$ $\left( \text{or, more universally, }h < \sqrt{\frac{\RE_1}{12\pi^2c_0\RE_2J G }}\right)$, then with the following constants:
\begin{itemize}
\item $C:= 32\pi^2\frac{1}{\RE_1}G_1^2\Gamma'\frac{e^{\gamma/2}}{e-1}$
\item $\gamma := \liminf\limits_{t\to \infty} \int\limits_t^{t+T} \frac{1}{2}\left[\mu-\frac{J^2}{\mu}|A\bu|^2\right] d\tau$,
\item  $\Gamma := \limsup\limits_{t\to \infty} \int\limits_t^{t+T} \max\left\{\frac{1}{2}\left[\mu-\frac{J^2}{\mu}|A\bu|^2\right],0\right\} d\tau$,
\item $\Gamma' :=e^{\Gamma+1+\gamma/2}$,
\end{itemize}
and for any $T \geq \frac{4\pi^2}{\RE_1}$, we obtain
\begin{align*}
 \|\bu(t)-\bv(t)\|^2 \leq \|\bu(0)-\bv(0)\|^2\Gamma ' e^{-\frac{\gamma}{2T}(t-t_0)} + C\cdot \RE_2(\RE_1^{-1}-\RE_2^{-1})^2.
\end{align*}
In particular,
\begin{align*}
\limsup\limits_{t\to\infty} \|\bu(t)-\bv(t)\| &\leq C \sqrt{\RE_2}|\RE_1^{-1}-\RE_2^{-1}|
= C\frac{|\RE_2-\RE_1|}{\RE_1\sqrt{\RE_2}}.
\end{align*}
\end{theorem}
\begin{proof}
We subtract \eqref{NSEmodified_mo} from \eqref{NSE_mo} to get
\begin{align*}
\bw_t + B(\bw,\bu)+B(\bv, \bw) = (\RE_2^{-1}-\RE_1^{-1})A\bu - \RE_2^{-1} A\bw - \mu P_\sigma( I_h(\bw)),
\end{align*}
with $\bw = \bu-\bv$ which, using the identity $B(\bw,\bu) - B(\bv,\bw) = B(\bu,\bw)+B(\bw,\bu)-B(\bw,\bw)$, can be simplified to 
\begin{align*}
&\frac{1}{2}\frac{d}{dt} \|\bw\|^2 + (B(\bu,\bw),A\bw) + (B(\bw,\bu),A\bw) - (B(\bw,\bw),A\bw) \\
&\qquad \qquad = (\RE_2^{-1}-\RE_1^{-1})(A\bu,A\bw)-\RE_2^{-1}|A\bw|^2 - \mu(P_\sigma(I_h(\bw)),A\bw). \notag
\end{align*}
Then, by \eqref{bilinear identity} and \eqref{lastbilinear}
\begin{align*}
\frac{1}{2}\frac{d}{dt} \|\bw\|^2 -(B(\bw,\bw),A\bu) =& (\RE_2^{-1}-\RE_1^{-1})(A\bu,A\bw)\\
& -\RE_2^{-1}|A\bw|^2 - \mu(P_\sigma(I_h(\bw)),A\bw). \notag
\end{align*}
Hence,
\begin{align*}
\frac{1}{2}\frac{d}{dt} \|\bw\|^2 + \RE^{-1}_2|A\bw|^2 = &(B(\bw,\bw), A\bu) \\
&+ (\RE_2^{-1}-\RE_1^{-1})(A\bu,A\bw)-\mu(P_\sigma (I_h(\bw)),A\bw).\notag
\end{align*}

By the Brezis-Gallouet inequality, 
\begin{align*}
|(B(\bw,\bw),A\bu)| \leq c\|\bw\|\left\{1+\log \frac{|A\bw|^2}{4\pi^2\|\bw\|^2}\right\}|A\bu|.
\end{align*}
Moreover, since $\mu c_0 h^2 \leq \RE_2^{-1}$ by assumption, we obtain
\begin{align*}
-\mu(P_\sigma(I_h(\bw)),A\bw) &= \mu(\bw-P_\sigma(I_h(\bw)),A\bw) - \mu\|\bw\|^2 \\
&\leq \mu|P_\sigma(\bw-I_h(\bw))||A\bw| - \mu\|\bw\|^2 \notag\\
&\leq \frac{\mu^2c_0h^2\RE_2}{2}\|\bw\|^2 + \frac{1}{2\RE_2}|A\bw|^2 - \mu\|\bw\|^2\notag\\
&\leq \frac{1}{2\RE_2}|A\bw|^2 - \frac{\mu}{2}\|\bw\|^2\notag.
\end{align*}
Using the Cauchy-Schwarz and Young's inequalities, we find
\begin{align*}
(\RE_2^{-1}-\RE_1^{-1})(A\bu,A\bw) &\leq |\RE_2^{-1}-\RE_1^{-1}||A\bu||A\bw| \\
&\leq \RE_2|\RE_2^{-1}-\RE_1^{-1}|^2|A\bu|^2 + \frac{1}{4\RE_2}|A\bw|^2.\notag
\end{align*}
Together, the above three inequalities imply that
\begin{align*}
\frac{1}{2}\frac{d}{dt} \|\bw\|^2 + \frac{3}{4\RE_2}|A\bw|^2 
&\leq c\|\bw\|^2\left(1+\log\frac{|A\bw|}{4\pi^2\|\bw\|^2}\right)|A\bu| + \RE_2(\RE_2^{-1}-\RE_1^{-1})^2|A\bu|^2 \notag\\
& \qquad + \frac{1}{4\RE_2}|A\bw|^2-\frac{\mu}{2}\|\bw\|^2\notag.
\end{align*}
Hence,
\begin{align*}
&\quad
\frac{d}{dt} \|\bw\|^2 + \frac{1}{\RE_2}|A\bw|^2  +\|\bw\|^2\left[\mu-2c|A\bu|\left(1+\log\frac{|A\bw|^2}{4\pi^2\|\bw\|^2}\right)\right] 
\\&\leq
2\RE_2(\RE_2^{-1}-\RE_1^{-1})^2|A\bu|^2.
\end{align*}

Let 
\begin{alignat*}{2}
\beta=\frac{\RE_2c|A\bu|}{2\pi^2} && \qquad \text{ and }\qquad   r = \frac{|A\bw|^2}{4\pi^2\|\bw\|^2}.
\end{alignat*}
Applying Lemma \ref{log inequality} (which is applicable since $r \geq 1$ by Poincar\'e's inequality), we obtain
\begin{align*}
\frac{-\RE_2c|A\bu|}{2\pi^2}\log\left(\frac{\RE_2c|A\bu|}{2\pi^2}\right) \leq \frac{|A\bw|^2}{4\pi^2\|\bw\|^2} - \frac{\RE_2c|A\bu|}{2\pi^2}\left(1+\log\frac{|A\bw|^2}{4\pi^2\|\bw\|^2}\right)
\end{align*}
which can be simplified to
\begin{align*}
-2c|A\bu|\log \left(\frac{\RE_2c|A\bu|}{2\pi^2}\right) \leq \frac{1}{\RE_2}\frac{|A\bw|^2}{\|\bw\|^2} - 2c|A\bu|\left(1+\log\frac{|A\bw|^2}{4\pi^2\|\bw\|^2}\right).
\end{align*}

This implies that 
\begin{align*}
&\frac{1}{2}\frac{d}{dt} \|\bw\|^2 + \|\bw\|^2\left[\mu-2c|A\bu|\log \frac{\RE_2c|A\bu|}{2\pi^2}\right] \leq 2\RE_2(\RE_2^{-1}-\RE_1^{-1})^2|A\bu|^2.\notag
\end{align*}

By \eqref{Au L2 bound}, 
\begin{align*}
2c\log \frac{\RE_2c|A\bu|}{2\pi^2} &\leq 2c\log\left(\frac{\RE_2c}{2\pi^2} \cdot \frac{\sqrt{c}4\pi^2}{\RE_1}(1+G)^2\right) \notag\\
&= 2c\log\left(\frac{2\RE_2c^{3/2}}{\RE_1}\right) + 4c\log(1+G).
\end{align*}
Let $J := 2c\log\left(\frac{2\RE_2c^{3/2}}{\RE_1}\right) + 4c\log(1+G)$.  Then
\begin{align*}
\frac{d}{dt}\|\bw\|^2 + \left[\mu-J|A\bu|\right]\|\bw\|^2 \leq 2\RE_2(\RE_2^{-1}-\RE_1^{-1})^2|A\bu|^2.
\end{align*}
Young's inequality implies
\begin{align*}
J|A\bu| \leq \frac{J^2}{2\mu}|A\bu|^2+\frac{\mu}{2}
\end{align*}
hence,
\begin{align*}
\frac{d}{dt}\|\bw\|^2 + \frac{1}{2}\left[\mu-\frac{J^2}{\mu}|A\bu|^2\right]\|\bw\|^2 \leq 2\RE_2(\RE_2^{-1}-\RE_1^{-1})^2|A\bu|^2.
\end{align*}

Next, we denote $\alpha(t) := \frac{1}{2}\left[\mu-\frac{J^2}{\mu}|A\bu|^2\right]$ and let $T := \frac{4\pi^2}{\RE_1}$.  Thanks to Theorem \ref{ubounds} and the assumption $\mu \geq 12\pi^2 \RE_1^{-1}JG_1$, it follows that
\begin{align*}
\gamma := \liminf\limits_{t\to \infty} \int\limits_t^{t+T} \alpha(\tau) d\tau &= \frac{\mu\RE_1}{8\pi^2} - \frac{J^2}{2\mu}
=\liminf_{t\to\infty}\int\limits_t^{t+T}|A\bu|^2 d\tau
\\&
\geq \frac{\mu\RE_1}{8\pi^2} - \frac{J^2}{2\mu}(\frac{16\pi^2G^2}{\RE_1}) 
> \frac{3}{2} JG - \frac{2}{3}JG 
= \frac{5}{6}JG > 0.
\end{align*}
Clearly, it follows that
\begin{align*}
\Gamma := \limsup\limits_{t\to\infty}\int\limits_t^{t+T} \alpha^-(\tau) d\tau < \infty,
\end{align*}
where $\alpha^-(t)$ is defined as in Lemma \ref{GronwallLemma}.

Choose $t_0$ sufficiently large so that Theorem \ref{ubounds} holds and the inequalities \eqref{alphaGamma} and \eqref{alphagamma} hold. Then, we can apply Lemma \ref{GronwallLemma} to obtain, for a.e. $t>t_0$,
\begin{align*}
 \|\bw(t)\|^2 &\leq \|\bw(t_0)\|^2\Gamma' e^{-\frac{\gamma}{2T}(t-t_0)}+\left(\sup\limits_{t\geq t_0}\int\limits_{t}^{t+T} 2\RE_2 |\RE_2^{-1}-\RE_1^{-1}|^2 |A\bu(\tau)|^2\; d\tau \right) \Gamma' \frac{e^{\gamma/2}}{e-1},
\end{align*}
where $\Gamma'$ as defined in Lemma \ref{GronwallLemma}.  
Taking the limit supremum as $t\rightarrow0$ establishes the result.
%
\end{proof}

Sensitivity for partial differential equations has been studied formally in many contexts (see, e.g., \cite{Anderson_Newman_Whitfield_Nielsen_1999_AIAA,Borggaard_Burns_1997,Brewer_1982_JMathAnal,Davis_Pahlevani_2013,Fernandez_Moubachir_2002_MMMAS,Gibson_Clark_1977,Hamby_1994_EnvMonAssess,Hyoungjin_Chongam_Rho_DongLee_1999_KSIAM,Kouhi_Houzeaux_Cucchietti_Vazquez_2016_AIAAConf,Neda_Pahlevani_Rebholz_Waters_2016,Pahlevani_2004,Pahlevani_2006,Rebholz_Zerfas_Zhao_2017_JMFM,Stanley_Stewart_2002,Vemuri_Raefsky_1979_IJSS}.  

\section{Sensitivity}\label{secSensitivity}
In this section, we analyze the sensitivity of $\bw$ to the Reynolds number by considering individually the sensitivity of $\bu$ and $\bv$ to the Reynolds number.  We wish to consider taking a derivative of equations \eqref{NSE_mo} and \eqref{NSEmodified_mo} with respect to the Reynolds number.  This has been done formally in many works on sensitivity (see, e.g., \cite{Anderson_Newman_Whitfield_Nielsen_1999_AIAA, Borggaard_Burns_1997, Breckling_Neda_Pahlevani_2018_CMA, Davis_Pahlevani_2013, Fernandez_Moubachir_2002_MMMAS, Hamby_1994_EnvMonAssess, Hyoungjin_Chongam_Rho_DongLee_1999_KSIAM, Kouhi_Houzeaux_Cucchietti_Vazquez_2016_AIAAConf,Pahlevani_2004, Pahlevani_2006, Vemuri_Raefsky_1979_IJSS}), yielding what are known as the \textit{sensitivity equations}.  However, to the best of our knowledge, a rigorous treatment has yet to appear in the literature.  Therefore, we provide a rigorous justification here of the existence and uniqueness of weak and strong solutions to the sensitivity equations in the case of zero initial data, which is the natural data for the sensitivity equation, as discussed below.  Moreover, we prove that these solutions can be realized as limits of difference quotients of Navier-Stokes solutions with respect to different Reynolds numbers.  Indeed, this is the method of our existence proofs, rather than using, e.g., Galerkin methods, fixed-point methods, etc.  Proofs using limits of difference quotients have appeared in the literature before, such as in standard proofs of elliptic regularity, the corresponding result for the Stokes equations, etc.  However, in the present context (i.e., the time-dependent sensitivity equations for 2D Navier-Stokes), we believe such a proof strategy is novel.

Working formally for a moment, we take the derivative of \eqref{NSE_mo} with respect to $\RE$, and denote (again, formally) $\widetilde{\bu} := \frac{d\bu_1}{d(\RE_1^{-1})}$ and $\widetilde{p} := \frac{dp_1}{d(\RE_1^{-1})}$, to obtain
\begin{subequations}\label{NSEpreSens1}
\begin{alignat}{2}
 \widetilde{\bu}_t + \widetilde{\bu}\cdot \nabla \bu_1 + \bu_1\cdot \nabla \widetilde{\bu} - &\RE_1^{-1} \triangle \widetilde{\bu} - \triangle{
 \bu_1} + \nabla \widetilde{p} = 0, \\
 \nabla \cdot \widetilde{\bu} &= 0.
\end{alignat}
\end{subequations}
These are known as the sensitivity equations for the Navier-Stokes equations.  
Similarly, denoting $\widetilde{\bv} := \frac{d\bv_1}{d(\RE_1^{-1})}$ and $\widetilde{q} := \frac{dq_1}{d(\RE_1^{-1})}$, we formally obtain
\begin{subequations}\label{NSEmodifiedSens1}
\begin{alignat}{2}
 \widetilde{\bv}_t + \widetilde{\bv}\cdot \nabla \bv_1 + \bv_1\cdot \nabla \widetilde{\bv} - &\RE_1^{-1} \triangle \widetilde{\bv} - \triangle{
 \bv_1} + \nabla \widetilde{q} = \mu I_h(\widetilde{\bu}-\widetilde{\bv}), \\
 \nabla \cdot \widetilde{\bv} &= 0,
\end{alignat}
\end{subequations}
Below, we prove some well-posedness results for these systems in the case of zero initial data.  We begin by defining what we mean by solutions.

\begin{remark}
We note that the sensitivity equations are a model for the evolution of the instantaneous change in a solution with respect to changes in the (inverse) Reynolds number.  Therefore, the natural initial condition to consider is the case of identically-zero initial data.  Indeed, if the initial data for the sensitivity equations is not identically zero, this would correspond to the case where the initial data for the Navier-Stokes equations depends on the viscosity, which is not typical of most mathematical treatments of the Navier-Stokes equations.  Thus, although we define weak solutions for general initial data, we only prove their existence for initial data which is identically zero.  Existence for general initial data can be proved using, e.g., Galerkin methods.  However, since our main focus is not on existence, but on showing the solutions can be realized as limits of a (sub)sequence of difference quotients, and moreover since the initial data is naturally taken to be zero in this setting, we use the difference quotient method instead.
\end{remark}

\begin{definition}\label{defweaksolnu} 
Let $T>0$.  Let $\bu\in L^2(0,T;V) \cap C_w(0,T;H)$ be a weak solution to \eqref{NSEpre} with initial data $\bu_0 \in V$ and forcing $\mathbf{f} \in L^\infty(0,\infty;V^*)$. A \textit{weak solution} of \eqref{NSEpreSens1} is an element $\widetilde{\bu} \in L^2(0,T;V) \cap C_w(0,T;H)$ satisfying $\frac{d\widetilde{\bu}}{dt} \in L^1_{\text{loc}}(0,T; V^*)$ and
 \begin{align}\label{NSE_sens_weak_form}
  \ip{\widetilde{\bu}_t}{\mathbf{\phi}} + \ip{B(\widetilde{\bu}, \bu)}{ \mathbf{\phi}} + \ip{B(\bu, \widetilde{\bu})}{ \mathbf{\phi}}+\RE_1^{-1}\ip{A\widetilde{\bu}}{\mathbf{\phi}} + \ip{A\bu}{\mathbf{\phi}} = 0
 \end{align}
 for a.e. $t\in[0,T]$, for all $\mathbf{\phi} \in V$, and initial data $\widetilde{\bu}_0 \in H$, satisfied in the sense of $C_w(0,T;H)$. 
 
 If, in addition, $\mathbf{f} \in L^\infty(0,\infty;H)$, $\bu_0 \in V$, $\widetilde{\bu}_0 \in V$, and $\bu\in L^2(0,T;V) \cap C_w(0,T;H)$ is a strong solution to \eqref{NSEpre}, then we define a \textit{strong solution} of \eqref{NSEpreSens1} to be a weak solution such that $\widetilde{\bu} \in L^2(0,T;\mathcal{D}(A)) \cap C^0([0,T]; V)$ and $\frac{d\widetilde{\bu}}{dt} \in L^2(0,T; H)$, satisfying \eqref{NSE_sens_weak_form} for a.e. $t \in [0,T]$ and for all $\phi \in H$.
\end{definition}

For the reasons discussed in Remark \ref{remark_no_weak_AOT_solutions} below, we only give a definition of strong solutions for the assimilation equations.

\begin{definition}\label{defweaksolnv} 
Let $T>0$.  Let $\bv$ be a strong solution to \eqref{NSEmodified} with initial data $\bv_0 \in V$ and forcing $\mathbf{f} \in L^\infty(0,\infty;H)$.
A \textit{strong solution} of \eqref{NSEmodifiedSens1} is an element $\widetilde{\bv} \in L^2(0,T;\mathcal{D}(A)) \cap C^0([0,T]; V)$ that satisfies 
 \begin{align*}
  \ip{\widetilde{\bv}_t}{\mathbf{\phi}} + \ip{B(\widetilde{\bv}, \bv)}{\mathbf{\phi}} + &\ip{B(\bv, \widetilde{\bv})}{\mathbf{\phi}}+\RE_1^{-1}\ip{A\widetilde{\bv}}{\mathbf{\phi}}\\ &+ \ip{A\bv}{\mathbf{\phi}} = \mu\ip{I_h(\widetilde{\bu}-\widetilde{\bv})}{\mathbf{\phi}}
 \end{align*}
 this equation for a.e. $t \in [0,T]$ and for all $\phi \in H$, where  $\frac{d\widetilde{\bv}}{dt} \in L^2(0,T; H)$ and initial data $\widetilde{\bv}_0 \in V$.
\end{definition}

Before we prove the existence and uniqueness of solutions with zero initial data to these equations, we first consider equations for the difference quotients.  Note that, since these are simple arithmetic operations on the Navier-Stokes equations, the manipulations can be performed rigorously, not just formally.  To this end, let $(\bu_1,p_1)$ be a solution to \eqref{NSEpre} with Reynolds number $\RE_1$ and $(\bu_2,p_2)$ be a solution to \eqref{NSEpre} with Reynolds number $\RE_2$ with the same initial data.  We take the difference of the two versions of \eqref{NSEpre}, each with Reynolds numbers $\RE_1$ and $\RE_2$.  We then divide by the difference in (inverse) Reynolds numbers, yielding the system
\begin{subequations}\label{NSEpreSens}
\begin{alignat}{2}
 \bD_t + \bu_2\cdot \nabla \bD + \bD \cdot \nabla \bu_1 - &\RE_2^{-1} \triangle \bD - \triangle \bu_1 + \nabla P = 0, \\
 \nabla \cdot \bD &= 0, \\
 \bD(\bx,0) &= 0,
\end{alignat}
\end{subequations}
where $\bD = \frac{\bu_1-\bu_2}{\RE_1^{-1}-\RE_2^{-1}}$ and $P := \frac{p_1-p_2}{\RE_1^{-1}-\RE_2^{-1}}$.  As defined, $\bD$ is a strong solution to \eqref{NSEpreSens}, and note that $\bu_1 = (\RE_1^{-1}-\RE_2^{-1}) \bD + \bu_2$.  Additionally, $\bD \in L^2(0,T; \mathcal{D}(A)) \cap C^0([0,T];V)$ and $\bD_t \in L^2(0,T; H)$. However, we need to establish that $\bD$ is the unique solution to \eqref{NSEpreSens}, which is the content of Lemma \ref{D_unique} below. 

\begin{lemma}\label{D_unique}
 Let $T>0$ be given, and let $\bu_1$, $\bu_2\in L^2(0,T; \mathcal{D}(A)) \cap C^0([0,T];V)$ be strong solutions to \eqref{NSEmodified}, with Reynolds numbers $\RE_1$ and $\RE_2$, respectively.  There exists one and only one solution $\bD$ to \eqref{NSEpreSens} that lies in $L^2(0,T; \mathcal{D}(A)) \cap C^0([0,T];V)$, i.e. for all $\phi \in H$, 
 \begin{align*}
(\bD_t,\phi)  + (B(\bD,\bu_1),\phi) + (B(\bu_2,\bD),\phi) + &\RE_2^{-1} (A\bD,\phi) +(A \bu_1,\phi) = 0,
 \end{align*}
 where $\bD_t \in L^2(0,T;H)$.
\end{lemma}

Next, we consider difference quotients for the assimilation system \eqref{NSEmodified}.  Let $(\bv_1,q_1)$ be the solution to \eqref{NSEmodified} with Reynolds number $\RE_1$ and $(\bv_2,q_2)$ be the solution to \eqref{NSEmodified} with Reynolds number $\RE_2$.  Subtracting the two equations and dividing by the difference in the (inverse) Reynolds numbers yields the system \eqref{NSEmodifiedSens},
\begin{subequations}\label{NSEmodifiedSens}
\begin{alignat}{2}
\bD'_t + \bD'\cdot \nabla \bv_1 + \bv_2\cdot \nabla \bD' - &\RE_2^{-1} \triangle \bD' - \triangle \bv_1 + \nabla Q = \mu I_h(\bD-\bD') \\
\nabla \cdot \bD' &= 0\\
\bD'(\bx,0) &= 0,
\end{alignat}
\end{subequations}
where $\bD' := \frac{\bv_1-\bv_2}{\RE_1^{-1}-\RE_2^{-1}}$ and $Q := \frac{q_1-q_2}{\RE_1^{-1}-\RE_2^{-1}}$.  As defined, $\bD'$ is a strong solution to \eqref{NSEmodifiedSens}, and note that $\bv_1 = (\RE_1^{-1}-\RE_2^{-1})\bD' + \bv_2$.  Additionally, $\bD' \in L^2(0,T; \mathcal{D}(A)) \cap C^0([0,T];V)$ and $\bD'_t \in L^2(0,T; H)$.

\begin{lemma}\label{D_prime_unique}
Let $T>0$ be given ,and let $\bv_1$, $\bv_2\in L^2(0,T; \mathcal{D}(A)) \cap C^0([0,T];V)$ be strong solutions to \eqref{NSEmodified}, with Reynolds numbers $\RE_1$ and $\RE_2$, respectively.  There exists a unique strong solution $\bD$ to \eqref{NSEmodifiedSens} that lies in $L^2(0,T; \mathcal{D}(A)) \cap C^0([0,T];V)$, in the sense that for all $\phi \in H$, 
 \begin{align*}
(\bD_t',\phi) + (B(\bv_2,\bD'),\phi) + (B(\bD',\nabla \bv_1),\phi) + &\RE_2^{-1} (A\bD',\phi) +(A \bv_1,\phi) = \mu(P_\sigma I_h(\bD-\bD'),\phi),
 \end{align*}
 where $\bD'_t \in L^2(0,T;H)$.
\end{lemma}

\begin{remark}\label{remark_no_weak_AOT_solutions}
The proofs of the above two lemmata are very similar; hence, we only present the proof of Lemma \ref{D_prime_unique}.  Moreover, we also note that in the case $\mu=0$, the proof of Lemma \ref{D_unique} holds \textit{mutatis mutandis} in the case where $\bu_1$, $\bu_2 \in C^0([0,T];H)\cap L^2(0,T; V)$ are only assumed to be weak solutions to the 2D Navier-Stokes equations, and then one obtains uniqueness of weak solutions to \eqref{NSEpreSens} in the class $C^0([0,T];H)\cap L^2(0,T;V)$.  However, in the case $\mu>0$, the notion of weak solutions for the assimilation equations \eqref{NSEmodified} has not been established in the literature for general interpolants $I_h$, and therefore we assume that the solutions $\bv_1$ and $\bv_2$ are strong solutions to \eqref{NSEmodified}, and prove the uniqueness of strong solutions to \eqref{NSEpreSens}.
\end{remark}

\begin{proof}
 Suppose there exist two solutions $\bD'_1$ and $\bD'_2$.  We consider the difference of the equations
\begin{align}
 \frac{d}{dt} \bD'_1 + B(\bD'_1, \bv_1) + B(\bv_2, \bD'_1) + \RE_2^{-1} A\bD'_1 + A\bv_1 = \mu P_\sigma I_h(\bD-\bD'_1)
\end{align}
and
\begin{align}
 \frac{d}{dt} \bD'_2 + B(\bD'_2, \bv_1) + B(\bv_2, \bD'_2) + \RE_2^{-1} A\bD'_2 + A\bv_1 = \mu P_\sigma I_h(\bD-\bD'_2)
\end{align}
which, defining $\bV := \bD'_1-\bD'_2$, yields
\begin{align}
 \bV_t + B(\bV,\bv_1) + B(\bv_2, \bV) + \RE_2^{-1} A\bV = -\mu P_\sigma I_h(\bV)
\end{align}
with $\bV(0) = 0$.  So, $\bV$ must be a solution to the above equation. Taking the inner product with $\bV$, 
\begin{align}
\frac{1}{2} \frac{d}{dt} |\bV|^2 + b(\bV,\bv_1,\bV) + \RE_2^{-1}\|\bV\|^2 = \ip{-\mu P_\sigma I_h(\bV)}{\bV}
\end{align}
which implies, applying the triangle inequality and Poisson's inequality to the interpolant term as in \cite{Azouani_Olson_Titi_2014},
\begin{align}
 &\quad
 \frac{1}{2} \frac{d}{dt} |\bV|^2 + \RE_2^{-1}\|\bV\|^2 
 \\&\leq\notag 
 c\|\bv_1\||\bV|\|\bV\| + \mu(\sqrt{c_0}h + \lambda_1^{-1})\|\bV\||\bV|
 \\&\leq\notag
 \frac{\mu^2(\sqrt{c_0}h + \lambda_1^{-1})^2}{\RE_2^{-1}} |\bV|^2 + \frac{\RE_2^{-1}}{4} \|\bV\|^2 + \frac{c^2}{2\RE_2^{-1}}\|\bv_1\|^2|\bV|^2 + \frac{\RE_2^{-1}}{2}\|\bV\|^2.
\end{align}
Thus, 
\begin{align}
 \frac{d}{dt} |\bV|^2 &\leq \Big( \frac{\mu^2(\sqrt{c_0}h + \lambda_1^{-1})^2}{\RE_2^{-1}} + \frac{c^2}{2\RE_2^{-1}}\|\bv_1\|^2\Big)|\bV|^2
\end{align}
and Gr{\"o}nwall's inequality implies
\begin{align}
 |\bV(T)|^2 \leq |\bV(0)|^2 \text{exp}\Big(\int_0^T  \frac{\mu^2(\sqrt{c_0}h + \lambda_1^{-1})^2}{\RE_2^{-1}} + \frac{c^2}{2\RE_2^{-1}}\|\bv\|^2 dt \Big).
\end{align}
But $\bV(0) = 0$, and thus $\|\bV\|_{L^\infty(0,T;H)} = 0$ implies that $\bV \equiv 0$.  Hence, solutions to \eqref{NSEmodifiedSens} are unique.
\end{proof}

Since systems \eqref{NSEpreSens} and \eqref{NSEmodifiedSens} have unique strong solutions for every $\RE_2^{-1}>0$, we want to show that, as $\RE_2 \to \RE_1$, the solutions to these equations converge to the unique strong solutions of the respective equations (in the sense of Definitions \ref{defweaksolnu} and \ref{defweaksolnv}) of the formal sensitivity equations \eqref{NSEpreSens1} and \eqref{NSEmodifiedSens1} with $0$ initial data.  We additionally prove that weak solutions exist for the sensitivity equations \eqref{NSEpreSens1} with $0$ initial data.

\begin{theorem}\label{ohhappyday1weak}
Let $\{(\RE_2^{-1})_n\}_{n \in \mathbb{N}}$ be a sequence such that $(\RE_2^{-1})_n \to \RE_1^{-1}$ as $n \to \infty$.  Let
\begin{itemize}
 \item  $\bu$ be a solution to \eqref{NSEpre} with Reynolds number $\RE_1^{-1}$, forcing $\mathbf{f} \in L^\infty(0,\infty;H)$, and initial data $\bu_0\equiv\mathbf{0}$;
 \item $\bu_2^n$ solve \eqref{NSEpre} with viscosity $(\RE_2^{-1})_n$, forcing $\mathbf{f} \in L^\infty(0,\infty;H)$, and initial data $\bu_0 \in V$;
 \item $\{\bD^n\}_{n \in \mathbb{N}}$ be a sequence of strong solutions to \eqref{NSEpreSens} with $\bD^n(0) = 0$.
 \end{itemize}
 Then there is a subsequence of $\{\bD^n\}_{n \in \mathbb{N}}$ that converges in $L^2(0,T;H)$ to a unique weak solution $\bD$ of \eqref{NSEpreSens1} with $0$ initial data for any $T>0$.
\end{theorem}
\begin{proof}
Let $T>0$ be given. Let $N$ sufficiently large such that for all $n>N$, $\{(\RE_2^{-1})_n\}_{n\in\mathbb{N}} \subset (\frac{\RE_1^{-1}}{2}, \frac{3\RE_1^{-1}}{2})$.  Then, we can follow the proof of strong solutions for \eqref{NSEpre} as in, e.g., \cite{Constantin_Foias_1988, Foias_Manley_Rosa_Temam_2001, Robinson_2001, Temam_2001_Th_Num}, to obtain bounds on $\{\bu_2^n\}$ for $n > N$ in the appropriate spaces that are independent of $(\RE_2^{-1})_n$:
\begin{align*}
\|\bu_2^n\|^2_{L^\infty(0,T;V)} &\leq \|\bu_2^n(0)\|^2 + \frac{ \|f\|^2_{L^2(0,T;H)}}{(\RE_2^{-1})_n} \\
&\leq  \|\bu_0\|^2 + \frac{ 2\|f\|^2_{L^2(0,T;H)}}{\RE_1^{-1}}
\end{align*}
and
\begin{align*}
\|\bu_2^n\|^2_{L^2(0,T;\mathcal{D}(A))} &\leq \frac{1}{(\RE_2^{-1})_n}\|\bu_2^n(0)\|^2 + \frac{ \|f\|^2_{L^2(0,T;H)}}{(\RE_2^{-1})_n^2} \\
&\leq \frac{2}{\RE_1^{-1}}\|\bu_0\|^2 + \frac{ 4\|f\|^2_{L^2(0,T;H)}}{(\RE_1^{-1})^2}.
\end{align*}
Note that $\|f\|^2_{L^2(0,T;H)} < \infty$ since all bounded functions are locally integrable.  Hence there is a subsequence that is relabeled $\bu_2^n \to \bu$ in $L^2(0,T;V)$ for some function $\bu$.  Continuing to follow the proof of strong solutions for \eqref{NSEpre} as in e.g. \cite{Constantin_Foias_1988, Foias_Manley_Rosa_Temam_2001, Robinson_2001, Temam_2001_Th_Num}, we note that $\frac{d \bu_2^n}{dt}$ is uniformly bounded in $n$ in $L^2(0,T;H)$.  Hence, we can find a subsequence which we relabel $\{\bu_2^n\}$ such that
\begin{align*}
 \frac{d \bu_2^n}{dt} \rightharpoonup \frac{d \bu}{dt} \qquad &\text{ in } L^2(0,T;H) \\
(\RE_2^{-1})_n A\bu_2^n \rightharpoonup \RE_1^{-1} A\bu \qquad &\text{ in } L^2(0,T;H) \\
B(\bu_2^n, \bu_2^n) \rightharpoonup B(\bu,\bu) \qquad &\text{ in } L^2(0,T;H).
\end{align*}
Indeed, $\bu$ satisfies \eqref{NSEpre} with corresponding Reynolds number $\RE_1^{-1}$ and thus, by uniqueness and the fact that $\bu_2^n \to \bu$ in $V$, $\bu_1 = \bu$.  Due to Poincar{\'e}'s inequality, we also obtain that $\bu_2^n \to \bu_1$ in $L^2(0,T;H)$.

Let $\bD^n$ be the strong solution to \eqref{NSEpreSens} with $\RE_1^{-1} = (\RE_2^{-1})_n$.  Taking the action of \eqref{NSEpreSens} on $\bD^n$ and using H{\"o}lder's, the bilinear inequalities, and Young's inequality twice, we obtain
\begin{align*}
 \frac{1}{2} \frac{d}{dt} |\bD^n|^2 + (\RE_2^{-1})_n \|\bD^n\|^2 &\leq 
 \frac{c^2}{(\RE_2^{-1})_n} \|\bu_1\|^2 |\bD^n|^2 + \frac{(\RE_2^{-1})_n}{4}\|\bD^n\|^2 \\
 &\phantom{=} + \frac{1}{2(\RE_2^{-1})_n}\|\bu_1\|^2 + \frac{(\RE_2^{-1})_n}{2}\|\bD^n\|^2,
\end{align*}
giving
\begin{align}\label{mainbound_Dn}
 \frac{1}{2}\frac{d}{dt}|\bD^n|^2 + \frac{(\RE_2^{-1})_n}{4}\|\bD^n\|^2 \leq  \frac{c^2}{(\RE_2^{-1})_n} \|\bu_1\|^2 |\bD^n|^2 + \frac{1}{2(\RE_2^{-1})_n}\|\bu_1\|^2.
\end{align}

Dropping the second term on the left hand side, we obtain
\begin{align*}
  \frac{1}{2}\frac{d}{dt}|\bD^n|^2 \leq  \frac{c^2}{(\RE_2^{-1})_n} \|\bu_1\|^2 |\bD^n|^2 + \frac{1}{2(\RE_2^{-1})_n}\|\bu_1\|^2.
\end{align*}

Taking the integral with respect to time on $[0,T]$ and applying Gr{\"o}nwall's inequality, then for a.e. $t \in [0,T]$,
\begin{align*}
 |\bD^n(t)|^2 &\leq \Big[\frac{1}{(\RE_2^{-1})_n}\int_0^T \|\bu_1\|^2 dt \Big] \text{exp}\Big(\int_0^T \frac{2c^2}{(\RE_2^{-1})_n} \|\bu_1\|^2 dt\Big) \\
 &\leq \Big[\frac{2}{\RE_1^{-1}}\int_0^T \|\bu_1\|^2 dt \Big] \text{exp}\Big(\int_0^T \frac{4c^2}{\RE_1^{-1}} \|\bu_1\|^2 dt\Big) =:K_1.
\end{align*}
Since $\bu_1 \in L^2(0,T; V)$, then $\bD^n$ is bounded above uniformly in $L^\infty(0,T;H)$.

Next, refraining from dropping the second term on the left hand side of \eqref{mainbound_Dn}, we estimate
\begin{align*}
 \frac{(\RE_2^{-1})_n}{4} \int_0^T \|\bD^n\|^2 dt &\leq \frac{c^2}{(\RE_2^{-1})_n} \int_0^T\|\bu_1\|^2 |\bD^n|^2dt + \frac{1}{2(\RE_2^{-1})_n}\int_0^T \|\bu_1\|^2dt \\
 &\leq  K_1\frac{c^2}{(\RE_2^{-1})_n} \int_0^T\|\bu_1\|^2 dt + \frac{1}{2(\RE_2^{-1})_n}\int_0^T \|\bu_1\|^2dt
\end{align*}
Rewriting, we obtain
\begin{align*}
 \int_0^T \|\bD^n\|^2 dt  &\leq K_1\frac{4c^2}{(\RE_2^{-1})^2_n} \int_0^T\|\bu_1\|^2 dt + \frac{2}{(\RE_2^{-1})^2_n}\int_0^T \|\bu_1\|^2dt \\
 &\leq K_1\frac{16c^2}{(\RE_1^{-1})^2} \int_0^T\|\bu_1\|^2 dt + \frac{8}{(\RE_1^{-1})^2}\int_0^T \|\bu_1\|^2dt
\end{align*}
Thus, $\bD^n$ is bounded above uniformly in $L^2(0,T;V)$ with respect to $n$.
Hence, by the Banach-Alaoglu Theorem, there exists a subsequence, relabeled as $(\bD^n)$, such that 
\begin{align}\label{BAconv_weak_u}
\bD^n \stackrel{*}{\rightharpoonup} \bD \text{ in } L^\infty(0,T;H)\quad \text{ and }\quad
\bD^n \rightharpoonup \bD \text{ in } L^2(0,T;V). 
\end{align}

Using \eqref{BAconv_weak_u}, note that all uniform bounds in $n$ on the terms in \eqref{NSEpreSens} in $L^2(0,T;V^*)$ are obtained in a similar manner to the proof of weak solutions for \eqref{NSEpre} except for the term $B(\bu_2^n,\bD^n)$.  However, by Lemma \ref{bilinear_unif_bd}, 
\begin{align*}
\|B(\bu_2^n,\bD^n)\|_{L^2(0,T;V^*)} &\leq k \|\bu_2^n\|_{L^\infty(0,T;H)}\|\bD^n\|_{L^\infty(0,T;H)}\|\bu_2^n\|_{L^2(0,T;V)}\|\bD^n\|_{L^2(0,T;V)},
\end{align*}
and due to the following bounds on $\bu_2^n$ (which can be found in \cite{Constantin_Foias_1988,Foias_Manley_Rosa_Temam_2001,Robinson_2001,Temam_2001_Th_Num}, etc.) and the fact that $(\RE_2^{-1})_n \in (\frac{\RE^{-1}}{2},\frac{3\RE^{-1}}{2})$,  
\begin{align*}
 \|\bu_2^n\|^2_{L^\infty(0,T;H)} &\leq |\bu^n_0|^2 + \frac{\|f\|_{L^\infty(0,T;H)}}{\lambda_1^2 (\RE_2^{-1})_n^2} \\
 &\leq |\bu_0|^2 + \frac{4\|f\|_{L^\infty(0,T;H)}}{\lambda_1^2 (\RE_1^{-1})^2}
\end{align*}
and
\begin{align*}
 \|\bu_2^n\|_{L^2(0,T;V)} &\leq \frac{1}{(\RE_2^{-1})_n}|\bu^n(0)|^2 + \frac{\|f\|_{L^\infty(0,T;H)}^2}{\lambda_1(\RE_2^{-1})_n^2} T \\
 &\leq \frac{2}{\RE_1^{-1}}|\bu_0|^2 + \frac{4\|f\|_{L^\infty(0,T;H)}^2}{\lambda_1(\RE_1^{-1})^2} T,
\end{align*}
 and thus $\|B(\bu_2^n,\bD^n)\|_{L^2(0,T;V^*)}$ is bounded above uniformly in $n$ independent of $(\RE_2^{-1})_n$.  Hence, independent of $(\RE_2^{-1})_n$, $d\bD^n/dt$ is bounded above uniformly in $n$ and by the Banach-Alaoglu Theorem a subsequence $\{\bD^n\}_{n\in\mathbb{N}}$ converges weakly to $d\bD/dt$ in $L^2(0,T;V^*)$.  Thus, by the Aubin Compactness Theorem, $\bD^n \to \bD$ strongly in $L^2(0,T;H)$.  Hence, weak continuity in $H$ follows due to the bounds on each of the terms above. Using these facts, we have weak-$*$ convergence in $L^2(0,T;V^*)$ of all but the bilinear terms in the standard sense.  Weak-$*$ convergence of the bilinear terms holds due to Lemma \ref{bilinear_wk_conv_v*},  yielding $B(\bD^n,\bu_1) \stackrel{*}{\rightharpoonup} B(\bD,\bu_1)$ in $L^2(0,T;V^*)$. Additionally since $\bu_2^n \to \bu_1$ strongly in $L^2(0,T;H)$, we can apply Lemma \ref{bilinear_wk_conv_v*} again to obtain that $B(\bu_2^n, \bD^n) \stackrel{*}{\rightharpoonup} B(\bu_1, \bD)$.  Thus, $\widetilde{\bu} := \bD$ satisfies
 \begin{align*}
  \widetilde{\bu}_t + B(\widetilde{\bu},\bu_1) + B(\bu_1, \widetilde{\bu}) + \RE_1^{-1} A\widetilde{\bu} + A\bu_1 = 0
 \end{align*}
in $L^2(0,T;V^*)$.
The initial condition is satisfied by construction.
 To prove uniqueness, suppose that there exist two weak solutions $\widetilde{\bu}_1$ and $\widetilde{\bu}_2$.  We consider the difference of the equations
\begin{align*}
 \frac{d}{dt} \widetilde{\bu}_1 + B(\widetilde{\bu}_1, \bu_1) + B(\bu_1,\widetilde{\bu}_1) + \RE_1^{-1} A\widetilde{\bu}_1 + A \bu_1  = 0
\end{align*}
and
\begin{align*}
 \frac{d}{dt} \widetilde{\bu}_2 + B(\widetilde{\bu}_2, \bu_1) + B(\bu_1, \widetilde{\bu}_2) + \RE_1^{-1} A\widetilde{\bu}_2 + A\bu_1 = 0,
\end{align*}
which, defining $\bU := \widetilde{\bu}_1 - \widetilde{\bu}_2$, yields
\begin{align*}
 \bU_t + B(\bU, \bu_1) + B(\bu_1, \bU) + \RE_1^{-1} A\bU = 0
\end{align*}
with $\bU(0) = 0$.  So, $\bU$ must be a weak solution to the above equation.  Taking the action on $\bU$ and applying the Lions-Magenes Lemma,

\begin{align*}
 \frac{1}{2}\frac{d}{dt} |\bU|^2 + \ip{B(\bU, \bu_1)}{\bU} + \RE_1^{-1}\|\bU\|^2 = 0 
\end{align*}
which implies
\begin{align*}
\frac{1}{2}\frac{d}{dt} |\bU|^2 +  \RE_1^{-1}\|\bU\|^2 &\leq c\|\bU\||\bU|\|\bu_1\| \\
&\leq \frac{c^2}{2\RE^{-1}} \|\bu_1\|^2|\bU|^2 + \frac{\RE_1^{-1}}{2} \|\bU\|^2.
\end{align*}
Dropping the second term, we obtain
\begin{align*}
 \frac{d}{dt} |\bU|^2 \leq \frac{c^2}{2\RE_1^{-1}} \|\bu_1\|^2|\bU|^2,
\end{align*}
and Gr{\"o}nwall's inequality implies that, for a.e. $0 \leq t \leq T$,
\begin{align*}
|\bU(t)|^2 \leq |\bU(0)|^2 \text{exp}\Big(\int_0^T \frac{c^2}{2\RE_1^{-1}} \|\bu_1\|^2 dt\Big). 
\end{align*}
Since we know the $\text{exp}\Big(\int_0^T \frac{c^2}{2\RE_1^{-1}} \|\bu_1\|^2 dt\Big)< \infty$ for all $T>0$ and $\bU(0) = 0$, we have that $\|\bU\|_{L^\infty(0,T; H)} = 0$, which implies that $\bU \equiv 0$.  Hence, weak solutions to \eqref{NSEpreSens} are unique.
\end{proof}

\begin{theorem}\label{ohhappyday1}
 Let $\{(\RE_2^{-1})_n\}_{n \in \mathbb{N}}$ be a sequence such that $(\RE_2^{-1})_n \to \RE_1^{-1}$ as $n \to \infty$.  Let
\begin{itemize}
 \item  $\bu$ be the solution to \eqref{NSEpre} with Reynolds number $\RE_1^{-1}$, forcing $\mathbf{f} \in L^\infty(0,\infty;H)$, and initial data $\bu_0$;
 \item $\bu_2^n$ solve \eqref{NSEpre} with viscosity $(\RE_2^{-1})_n$, forcing $\mathbf{f} \in L^\infty(0,\infty;H)$, and initial data $\bu_0 \in V$
 \item $\{\bD^n\}_{n \in \mathbb{N}}$ be a sequence of strong solutions to \eqref{NSEpreSens} with $\bD^n(0) = 0$.
 \end{itemize}
 Then there is a subsequence of $\{\bD^n\}_{n \in \mathbb{N}}$ that converges in $L^2(0,T;V)$ to a unique strong solution $\bD$ of \eqref{NSEpreSens1} with $0$ initial data.
\end{theorem}
\begin{proof}

Let $T>0$ be given, and let $N>0$ be large enough that $n>N$ implies $\{(\RE_2^{-1})_n\} \subset (\frac{\RE_1^{-1}}{2}, \frac{3\RE_1^{-1}}{2})$.  Then by the argument in Theorem \ref{ohhappyday1weak}, we can obtain a subsequence which we relabel $\{\bu_2^n\}$ such that $\bu_2^n \to \bu_1$ in $L^2(0,T;V)$.

Consider $\bD^n$ to be the strong solution to \eqref{NSEpreSens} with Reynolds number $(\RE_2^{-1})_n$.  Taking a justified inner product of \eqref{NSEpreSens} with $A\bD^n$,

\begin{align*}
\frac{1}{2} \frac{d}{dt} \|\bD^n\|^2 + (\RE_2^{-1})_n |A\bD^n|^2 &= -(B(\bD^n,\bu_1),A\bD^n) - (B(\bu_2^n, \bD^n), A\bD^n) \\
&\phantom{=}- (A\bu_1,A \bD^n).
\end{align*}
Applying Young's inequality, we obtain
\begin{align*}
\frac{1}{2} \frac{d}{dt} \|\bD^n\|^2 + \frac{(\RE_2^{-1})_n}{2} |A\bD^n|^2 &\leq -(B(\bD^n,\bu_1),A\bD^n) - (B(\bu_2^n, \bD^n), A\bD^n) \\
&\phantom{=} + \frac{1}{2(\RE_2^{-1})_n}|A\bu_1|^2.
\end{align*}
Applying \eqref{BINsimple} to the second bilinear term,
\begin{align*}
\frac{1}{2} \frac{d}{dt} \|\bD^n\|^2 + \frac{(\RE_2^{-1})_n}{2} |A\bD^n|^2 &\leq -(B(\bD^n,\bu_1),A\bD^n) + \|\bu_2^n\|_{L^\infty(\Omega)}\|\bD^n\||A\bD^n| \\
&\phantom{=} + \frac{1}{2(\RE_2^{-1})_n}|A\bu_1|^2 \\
&\leq \frac{2k^2}{(\RE_2^{-1})_n} |\bu_2^n||A\bu_2^n|\|\bD^n\|^2 + \frac{(\RE_2^{-1})_n}{8} |A\bD^n|^2 \\
&\phantom{=} -(B(\bD^n,\bu_1),A\bD^n)  + \frac{1}{2(\RE_2^{-1})_n}|A\bu_1|^2
\end{align*}
and applying \eqref{BIN2} to the first bilinear term,
\begin{align*}
\frac{1}{2} \frac{d}{dt} \|\bD^n\|^2 + \frac{3(\RE_2^{-1})_n}{8} |A\bD^n|^2 &\leq \frac{2k^2}{(\RE_2^{-1})_n}|\bu_2^n||A\bu_2^n|\|\bD^n\|^2 \\
&\phantom{=} + c|\bD^n|^{1/2}\|\bD^n\|^{1/2}\|\bu_1\|^{1/2}|A\bu_1|^{1/2}|A\bD^n| + \frac{1}{2(\RE_2^{-1})_n}|A\bu_1|^2 \\
&\leq \frac{2k^2}{(\RE_2^{-1})_n}|\bu_2^n||A\bu_2^n|\|\bD^n\|^2  + \frac{2c^2}{\lambda_1(\RE_2^{-1})_n} \|\bD^n\|^2\|\bu_1\||A\bu_1|\\
&\phantom{=}+ \frac{(\RE_2^{-1})_n}{8}|A\bD^n|^2  + \frac{1}{2(\RE_2^{-1})_n}|A\bu_1|^2
\end{align*}
which can be rewritten as 
\begin{align*}
 \frac{d}{dt} \|\bD^n\|^2 + \frac{(\RE_2^{-1})_n}{2} |A\bD^n|^2 &\leq \Big(\frac{4k^2}{(\RE_2^{-1})_n}|\bu_2^n||A\bu_2^n| +  \frac{4c^2}{\lambda_1(\RE_2^{-1})_n}\|\bu_1\||A\bu_1| \Big) \|\bD^n\|^2 \\
 &\phantom{=} + \frac{1}{(\RE_2^{-1})_n}|A\bu_1|^2.
\end{align*}
Integrating on both sides in time, with $0 \leq t \leq T$,
\begin{align*}
\|\bD^n(t)\|^2 + &\frac{(\RE_2^{-1})_n}{2} \int_0^t |A\bD^n|^2 ds \leq  \frac{1}{(\RE_2^{-1})_n}\int_0^t|A\bu_1(s)|^2 ds \\
 &\phantom{=} + \int_0^t\Big(\frac{4k^2}{(\RE_2^{-1})_n} |\bu_2^n(s)||A\bu_2^n(s)| +  \frac{4c^2}{\lambda_1(\RE_2^{-1})_n}\|\bu_1(s)\||A\bu_1(s)| \Big) \|\bD^n(s)\|^2 ds 
\end{align*}
Dropping the second term on the left hand side, we apply Gr{\"o}nwall's inequality to obtain
\begin{align*}
\|\bD^n(t)\|^2 &\leq \alpha_n(t) \;\text{exp} \Big(\int_0^t\frac{4k^2}{(\RE_2^{-1})_n} |\bu_2^n(s)||A\bu_2^n(s)| +  \frac{4c^2}{\lambda_1(\RE_2^{-1})_n}\|\bu_1(s)\||A\bu_1(s)| ds\Big) \\
&\leq \alpha(t)\;\text{exp} \Big(\int_0^t\frac{8k^2}{\RE_1^{-1}} |\bu_2^n(s)||A\bu_2^n(s)| +  \frac{8c^2}{\lambda_1\RE_1^{-1}}\|\bu_1(s)\||A\bu_1(s)| ds\Big).
\end{align*}
where $\alpha_n(t) := \frac{2}{(\RE_2^{-1})_n}\int_0^t|A\bu_1(s)|^2 ds \leq \alpha(t) := \frac{4}{\RE_1^{-1}}\int_0^t|A\bu_1(s)|^2 ds$.  Since 
\[ \int_0^T |A\bu_2^n|^2 ds \leq \|\bu_0\|^2 + \frac{\|f\|_{L^2(0,T; H)}}{(\RE_2^{-1})_n}\]
as proven in, e.g., \cite{Constantin_Foias_1988, Robinson_2001, Foias_Manley_Rosa_Temam_2001, Temam_2001_Th_Num}, then
\begin{align*}
 \sup\limits_{t\in [0,T]} \|\bD^n(t)\|^2 &\leq \alpha(T) \frac{8k^2}{\lambda_1^2(\RE_1^{-1})} \int_0^T |A\bu_2^n|^2 ds + \frac{8c^2}{\lambda_1(\RE_1^{-1})}\|\bu_1(s)\||A\bu_1(s)| ds\\
 &\leq \alpha(T) \frac{8k^2}{\lambda_1^2(\RE_1^{-1})}\Big[ \|\bu_0\|^2 + \frac{\|f\|_{L^2(0,T; H)}}{(\RE_2^{-1})_n}\Big] \\
 &\phantom{=} + \alpha(T) \int_0^T \frac{8c^2}{\lambda_1(\RE_1^{-1})}\|\bu_1(s)\||A\bu_1(s)| ds\\
 &\leq \alpha(T) \frac{8k^2}{\lambda_1^2(\RE_1^{-1})}\Big[ \|\bu_0\|^2 + \frac{2\|f\|_{L^2(0,T; H)}}{(\RE_1^{-1})}\Big] \\
 &\phantom{=} + \alpha(T) \int_0^T \frac{8c^2}{\lambda_1(\RE_1^{-1})}\|\bu_1(s)\||A\bu_1(s)| ds\\
\end{align*}

This implies that $\bD^n \in L^\infty(0,T; V)$ and $\{\bD^n\}$ is uniformly bounded in this space.

Additionally, considering again the inequality
\begin{align*}
\|\bD^n(t)\|^2 + &\frac{(\RE_2^{-1})_n}{2} \int_0^t |A\bD^n|^2 ds \leq \frac{1}{(\RE_2^{-1})_n}\int_0^t|A\bu_1(s)|^2 ds \\
 &\phantom{=} + \int_0^t\Big(\frac{4k^2}{(\RE_2^{-1})_n} |\bu_2^n(s)||A\bu_2^n(s)| +  \frac{4c^2}{\lambda_1(\RE_2^{-1})_n}\|\bu_1(s)\||A\bu(s)| \Big) \|\bD^n(s)\|^2 ds.
\end{align*}
we set $t = T$, drop the first term on the left hand side, and bound the Reynolds number above to obtain
\begin{align*}
\int_0^T |A\bD^n|^2 ds &\leq \frac{8}{(\RE_1^{-1})^2}\left(\int_0^T|A\bu_1(s)|^2 ds \right)\\
 &\phantom{=} + \int_0^T\Big(\frac{32k^2}{\lambda_1(\RE_1^{-1})^2} |A\bu_2^n(s)|^2 +  \frac{32c^2}{\lambda_1(\RE_1^{-1})^2}\|\bu_1(s)\||A\bu(s)| \Big) \|\bD^n(s)\|^2 ds
\end{align*}
By the fact that $\{\|\bu_2^n\|_{L^2(0,T;\mathcal{D}(A))}\}$ is bounded above in $n$ as demonstrated in Theorem \ref{ohhappyday1weak} and the result that $\{\|\bD^n\|_{L^\infty(0,T;V)}\}$ is bounded above uniformly in $n$, we also have that $\{\|\bD^n\|_{L^2(0,T; \mathcal{D}(A))}\}$ is bounded above uniformly in $n$.  Since $\{\bD^n\}$ is bounded above uniformly in $n$ in both $L^\infty(0,T;V)$ and $L^2(0,T;\mathcal{D}(A))$, then we can conclude that there exists a subsequence, which we relabel as $\{\bD^n\}$, such that 
\begin{equation}\label{BAconv_str_u}
\bD^n \stackrel{*}{\rightharpoonup} \bD \text{ in } L^\infty(0,T; V) \text{ and } \bD^n \rightharpoonup \bD \text{ in } L^2(0,T;\mathcal{D}(A)). 
\end{equation}
Using \eqref{BAconv_str_u}, note that all uniform bounds in $n$ on the terms in \eqref{NSEpreSens} in $L^2(0,T;H)$ are obtained in a similar manner to the proof of strong solutions for the \eqref{NSEpre} and are independent of $(\RE_2^{-1})_n$ except for the bilinear terms.  The bilinear terms are bounded uniformly in $L^2(0,T;H)$ with respect to $n$, due to Lemma \ref{bilinear_unif_bd}.  Hence, $\frac{d\bD^n}{dt}$ is bounded above uniformly in $n$ in $L^2(0,T;H)$.  Thus, as in, e.g., \cite{Robinson_2001, Constantin_Foias_1988, Foias_Manley_Rosa_Temam_2001, Temam_2001_Th_Num}, 
\[\frac{d\bD^n}{dt} \rightharpoonup \frac{d\bD}{dt} \text{ in } L^2(0,T;H).\]
Hence, by the Aubin Compactness Theorem, $\bD^n \to \bD$ strongly in $L^2(0,T;V)$.  As in, e.g., \cite{Robinson_2001, Temam_2001_Th_Num, Foias_Manley_Rosa_Temam_2001, Constantin_Foias_1988}, $\bD \in C^0(0,T;V)$. Using these facts, we have weak convergence in $L^2(0,T;H)$ for all except the bilinear terms in the standard sense.  Weak convergence of the bilinear terms holds due to Lemma \ref{bilinear_wk_conv}.  Hence, $\widetilde{\bu} := \bD$ satisfies
 \begin{align*}
  \widetilde{\bu}_t + B(\widetilde{\bu},\bu_1) + B(\bu_1, \widetilde{\bu}) + \RE_1^{-1} A\widetilde{\bu} + A\bu_1 = 0
 \end{align*}
in $L^2(0,T;H)$.

The initial condition is also satisfied by construction.
Uniqueness holds due to the results in Theorem \ref{ohhappyday1weak}.
\end{proof}

\begin{theorem}\label{ohhappyday2}
Let $\{(\RE_2^{-1})_n\}_{n \in \mathbb{N}}$ be a sequence such that $(\RE_2^{-1})_n \to \RE_1^{-1}$ as $n \to \infty$.  Choose $\mu$ and $h$ such that $4\mu c_0 h^2 \leq (\RE_2^{-1})_n \leq \frac{3\RE_1^{-1}}{2}$.  Let
\begin{itemize}
 \item  $\bv$ be the solution to \eqref{NSEmodified} with Reynolds number $\RE_1^{-1}$, forcing $\mathbf{f} \in L^\infty(0,\infty;H)$, and initial data $\bv_0$;
 \item $\bv_2^n$ solve \eqref{NSEmodified} with viscosity $(\RE_2^{-1})_n$, forcing $\mathbf{f} \in L^\infty(0,\infty;H)$, and initial data $\bv_0 \in V$;
 \item $\{\bD^{n'}\}_{n \in \mathbb{N}}$ be a sequence of strong solutions to \eqref{NSEmodifiedSens} with $\bD^{n'}(0) = 0$.
 \end{itemize}
 Then there is a subsequence of $\{\bD^{n'}\}_{n \in \mathbb{N}}$ that converges in $L^2(0,T;V)$ to a unique solution $\bD'$ of \eqref{NSEmodifiedSens1} with $0$ initial data.
\end{theorem}
\begin{proof}

 Let $T > 0$. Note that since $\{(\RE_2^{-1})_n\} \subset (\frac{\RE_1^{-1}}{2}, \frac{3\RE_1^{-1}}{2})$  for $n > N$ for some sufficiently large $N$, we can follow the proof of strong solutions for \eqref{NSEmodified} in \cite{Azouani_Olson_Titi_2014} to obtain bounds on $\{\bv_2^n\}_{n>N}$ in the appropriate spaces that are independent of $(\RE_2^{-1})_n$.  First, we note that \cite{Azouani_Olson_Titi_2014} quickly proves $|f+\mu P_\sigma I_h(\bu_2^n)|^2 \leq M_n$ since $|P_\sigma I_h(\bu_2^n)|^2 \leq |\bu_2^n|^2$.  However, since $\bu_2^n$ is bounded above uniformly in $n$ (see the proof of Theorem \ref{ohhappyday1weak}), we have that $|f+\mu P_\sigma I_h(\bu_2^n)|^2 \leq m$ for some $m$ independent of $n$.  Thus, we have the following bounds from \cite{Azouani_Olson_Titi_2014} bounded above uniformly in $n$:
\begin{align}\label{vbound_2_inf}
\|\bv_2^n\|_{L^\infty(0,T;H)}^2 &\leq |\bv^n_2(0)|^2 + \frac{m}{\mu (\RE_2^{-1})_n \lambda_1} \\
   &\leq |\bv_0|^2 + \frac{2m}{\mu \RE_1^{-1} \lambda_1} \notag,
\end{align}
\begin{align}\label{vbound_2_2}
\|\bv_2^n\|_{L^2(0,T;V)}^2 &\leq \frac{1}{(\RE_2^{-1})_n}|\bv^n_2(0)|^2 + \frac{T}{\mu(\RE_2^{-1})_n} m \\
 &\leq \frac{2}{\RE_1^{-1}}|\bv_0|^2 + \frac{2T}{\mu\RE_1^{-1}} m, \notag
\end{align}
\begin{align}\label{vbound_h1_inf}
 \|\bv_2^n\|_{L^\infty(0,T;V)}^2 &\leq \frac{1}{\psi(T)} \Big[ \|\bv^n_2(0)\|^2 + \frac{4T}{(\RE_2^{-1})_n}m\Big] \\
 &\leq \frac{1}{\overline{\psi(T)}} \Big[ \|\bv_0\|^2 + \frac{8T}{\RE_1^{-1}}m\Big]\notag
\end{align}
where 
\begin{align*}
\frac{1}{\psi(T)} &= \text{exp}\Big\{\frac{c}{(\RE_2^{-1})^3_n}\int_0^T|\bv_2^n|^2\|\bv_2^n\|^2 ds\Big\} \\
&\leq \frac{1}{\overline{\psi(T)}} = \text{exp}\Big\{\frac{8c}{(\RE_1^{-1})^3}\int_0^T|\bv_2^n|^2\|\bv_2^n\|^2 ds\Big\},
\end{align*}
which is bounded above uniformly in $n$ due to \eqref{vbound_2_inf} and \eqref{vbound_2_2}, and
\begin{align*}
\|\bv_2^n\|_{L^2(0,T;\mathcal{D}(A))}^2 &\leq \frac{1}{(\RE_2^{-1})_n}\|\bv^n_2(0)\|^2 + \frac{c}{(\RE_2^{-1})^3_n} \int_0^T (|\bv_2^n|^2\|\bv_2^n\|^4 + \frac{4}{(\RE_2^{-1})_n}|f + P_\sigma I_h(\bu_2^n)|^2) ds \\
 &\leq \frac{2}{(\RE_1)^{-1}}\|\bv_0\|^2 + \frac{8c}{(\RE_1^{-1})^3} \int_0^T |\bv_2^n|^2\|\bv_2^n\|^4 ds + \frac{8T}{\RE_1^{-1}}m,
\end{align*}
which is bounded above uniformly in $n$ due to \eqref{vbound_2_inf}, \eqref{vbound_2_2}, \eqref{vbound_h1_inf}.  Hence, we will obtain a subsequence that is relabeled $\bv_2^n \to \bv$ in $L^2(0,T;V)$ for some function $\bv$.  Indeed, we see that by identical arguments presented in Theorem \ref{ohhappyday1weak}, $\bv = \bv_1$.  Also due to Poincar{\'e}'s inequality, we obtain that $\bv_2^n \to \bv_1$ in $L^2(0,T;H)$.

Let $\{\bD^{n'}\}_{n \in \mathbb{N}}$ be a sequence of solutions to \eqref{NSEmodifiedSens}.  We consider the Leray projection of \eqref{NSEmodifiedSens}:
\begin{align*}
 \frac{d}{dt}\bD^{n'} + B(\bD^{n'},\bv_1) + B(\bv_2^n, \bD^{n'}) + (\RE_2^{-1})_nA \bD^{n'} + A \bv_1 = \mu P_\sigma I_h(\bD^n-\bD^{n'}).
\end{align*}
The existence proof for \eqref{NSEmodifiedSens1} closely follows the proof of Theorem \ref{ohhappyday1}, with some modifications on the bounds of $\bD^{n'}$ which we show below.  
Taking the inner product with $A\bD^{n'}$ and proceeding as in the proof of Theorem \ref{ohhappyday1}, we obtain
\begin{align}\label{dnprime_inequality}
 \frac{1}{2}\frac{d}{dt} \|\bD^{n'}\|^2 + \frac{(\RE_2^{-1})_n}{4} |A\bD^{n'}|^2 &\leq \Big(\frac{2k^2}{(\RE_2^{-1})_n}|\bv_2^n||A\bv_2^n| +  \frac{2c^2}{\lambda_1(\RE_2^{-1})_n}\|\bv_1\||A\bv_1| \Big) \|\bD^{n'}\|^2 
 \\\notag
 &\phantom{=} + \frac{1}{2(\RE_2^{-1})_n}|A\bv_1|^2 + \mu (I_h(\bD^n - \bD^{n'}),A\bD^{n'}).
\end{align}
We slightly modify the inequalities obtained in \cite{Azouani_Olson_Titi_2014} for the interpolant term,
\begin{align*}
\mu|(I_h(\bD^{n'}), A\bD^{n'})| &\leq \frac{4\mu^2}{(\RE_2^{-1})_n}|\bD^{n'}-I_h(\bD^{n'})|^2 + \frac{(\RE_2^{-1})_n}{16}|A\bD^{n'}|^2 - \mu\|\bD^{n'}\|^2 \\
&\leq \frac{4\mu^2c_0h^2}{(\RE_2^{-1})_n}\|\bD^{n'}\|^2 + \frac{(\RE_2^{-1})_n}{16}|A\bD^{n'}|^2 - \mu\|\bD^{n'}\|^2 \\
&\leq \frac{(\RE_2^{-1})_n}{16} |A\bD^{n'}|^2.
\end{align*}
Also,
\begin{align*}
\mu|(I_h(\bD^n), A\bD^{n'})| \leq \frac{4\mu^2}{(\RE_2^{-1})_n} |\bD^n|^2 + \frac{(\RE_2^{-1})_n}{16}|A\bD^{n'}|^2.
\end{align*}
Using these inequalities in \eqref{dnprime_inequality}:
\begin{align*}
 \frac{1}{2}\frac{d}{dt} \|\bD^{n'}\|^2 + \frac{(\RE_2^{-1})_n}{8} |A\bD^{n'}|^2 &\leq \Big(\frac{2k^2}{(\RE_2^{-1})_n}|\bv_2^n||A\bv_2^n| +  \frac{2c^2}{\lambda_1(\RE_2^{-1})_n}\|\bv_1\||A\bv_1| \Big) \|\bD^{n'}\|^2 \\
 &\phantom{=} + \frac{1}{2(\RE_2^{-1})_n}|A\bv_1|^2 + \frac{4}{(\RE_2^{-1})_n}|\bD^n|^2  \\
 &\leq \Big(\frac{2k^2}{(\RE_2^{-1})_n}|\bv_2^n||A\bv_2^n| +  \frac{2c^2}{\lambda_1(\RE_2^{-1})_n}\|\bv_1\||A\bv_1| \Big) \|\bD^{n'}\|^2 \\
 &\phantom{=} + \frac{1}{2(\RE_2^{-1})_n}|A\bv_1|^2 + \frac{4}{\lambda_1^2(\RE_2^{-1})_n}|A\bD^n|^2.
\end{align*}
Following identical arguments as in Theorem \ref{ohhappyday1} with \[\alpha_n(t) := \frac{1}{2(\RE_2^{-1})_n}|A\bv_1|^2 + \frac{4}{\lambda_1^2(\RE_2^{-1})_n}|A\bD^n|^2 \leq \alpha(t) := \frac{1}{\RE_1^{-1}}|A\bv_1|^2 + \frac{8}{\lambda_1^2\RE_1^{-1}}|A\bD^n|^2,\] along with the fact that $P_\sigma I_h(\bD^n-\bD^{n'})$ is bounded uniformly in $n$ in $L^2(0,T;H)$, we obtain a subsequence relabeled $\bD^{n'} \to \bD'$ in $L^2(0,T;V)$. Indeed, let $\phi \in L^2(0,T;H)$; then
\begin{align*}
\int_0^T (P_\sigma I_h(\bD^n-\bD^{n'}) &- P_\sigma I_h(\bD-\bD'),\phi) ds \leq \int_0^T |I_h(\bD^n-\bD^{n'}) - I_h(\bD-\bD')||\phi| ds \\
&\leq \int_0^T |I_h(\bD^n-\bD) - I_h(\bD^{n'}-\bD')||\phi| ds \\
&\leq \int_0^T |[(\bD^n-\bD) - (\bD^{n'}-\bD')] - I_h((\bD^n-\bD)-(\bD^{n'}-\bD'))||\phi| ds \\
&\phantom{=}+ \int_0^T |(\bD^n-\bD) - (\bD^{n'}-\bD')||\phi| ds \\
&\leq \sqrt{c_0}h \int_0^T \|(\bD^n-\bD) - (\bD^{n'}-\bD')\||\phi| ds \\
&\phantom{=}+ \frac{1}{\lambda_1^{1/2}} \int_0^T \|(\bD^n-\bD) - (\bD^{n'}-\bD')\||\phi| ds \\
&\leq (\sqrt{c_0}h \|\bD^n-\bD\|_{L^2(0,T;V)}\|\phi\|_{L^2(0,T;H)} \\
&\phantom{=} + \frac{1}{\lambda_1^{1/2}} \|\bD^{n'}-\bD'\|_{L^2(0,T;V)}\|\phi\|_{L^2(0,T;H)}).
\end{align*}

Additionally, since we now have that $\bD^{n'} \to \bD$ in $L^2(0,T;V)$, then $P_\sigma I_h(\bD^n-\bD^{n'}) \rightharpoonup P_\sigma I_h(\bD-\bD')$ in $L^2(0,T;H)$ and we conclude $\bD'$ is a strong solution in the sense of Definition \ref{defweaksolnv}. 

To show that the solutions are unique, we consider the 
difference of the equations
\begin{align*}
 \frac{d}{dt} \widetilde{\bv}_1 + B(\widetilde{\bv}_1, \bv_1) + B(\bv_1, \widetilde{\bv}_1) + \RE^{-1} A\widetilde{\bv}_1 + A\bv = \mu P_\sigma I_h(\widetilde{\bu}-\widetilde{\bv}_1)
\end{align*}
and
\begin{align*}
 \frac{d}{dt} \widetilde{\bv}_2 + B(\widetilde{\bv}_2, \bv_1) + B(\bv_1, \widetilde{\bv}_2) + \RE^{-1} A\widetilde{\bv}_2 + A\bv_1 = \mu P_\sigma I_h(\widetilde{\bu}-\widetilde{\bv}_2)
\end{align*}
which, defining $\bV := \widetilde{\bv}_1-\widetilde{\bv}_2$, yields
\begin{align*}
 \bV_t + B(\bV,\bv_1) + B(\bv_1, \bV) + \RE^{-1} A\bV = -\mu P_\sigma I_h(\bV)
\end{align*}
with $\bV(0) = 0$.  So, $\bV$ must be a solution to the above equation. Taking the action on $\bV$ and applying the Lions-Magenes Lemma, 
\begin{align*}
\frac{1}{2} \frac{d}{dt} |\bV|^2 + \ip{B(\bV,\bv_1)}{\bV} + \RE^{-1}\|\bV\|^2 = \ip{-\mu P_\sigma I_h(\bV)}{\bV}
\end{align*}
which implies that
\begin{align*}
 \frac{1}{2} \frac{d}{dt} |\bV|^2 + \RE^{-1}\|\bV\|^2 &\leq c\|\bv_1\||\bV|\|\bV\| + \mu(\sqrt{c_0}h + \lambda_1^{-1})\|\bV\||\bV| \\
 &\leq \frac{\mu^2(\sqrt{c_0}h + \lambda_1^{-1})^2}{\RE^{-1}} |\bV|^2 + \frac{\RE^{-1}}{4} \|\bV\|^2 \\
 &\phantom{=} + \frac{c^2}{2\RE^{-1}}\|\bv_1\|^2|\bV|^2 + \frac{\RE^{-1}}{2}\|\bV\|^2.
\end{align*}
Thus, 
\begin{align*}
 \frac{d}{dt} |\bV|^2 &\leq \Big( \frac{\mu^2(\sqrt{c_0}h + \lambda_1^{-1})^2}{\RE^{-1}} + \frac{c^2}{2\RE^{-1}}\|\bv_1\|^2\Big)|\bV|^2
\end{align*}
and Gr{\"o}nwall's inequality implies, for a.e. $0 \leq t \leq T$,
\begin{align*}
 |\bV(t)|^2 \leq |\bV(0)|^2 \text{exp}\Big(\int_0^T  \frac{\mu^2(\sqrt{c_0}h + \lambda_1^{-1})^2}{\RE^{-1}} + \frac{c^2}{2\RE^{-1}}\|\bv_1\|^2 dt \Big).
\end{align*}
But $\bV(0) = 0$, and thus $\|\bV\|_{L^\infty(0,T;H)} = 0$ implies that $\bV \equiv 0$.  Hence, solutions to \eqref{NSEmodifiedSens} are unique.
\end{proof}


\section{Numerical Results}\label{secNumerics}
In the previous sections, we showed that the data assimilation algorithm \eqref{NSEmodified} can still perform well even when there is error in the viscosity parameter, provided that $\mu$ is large and $h$ is small. However, for large values of the Grashof number satisfying the requirements of the rigorous estimates would require prohibitively small values of $h$. Fortunately, in practice the requirements on $h$ and $\mu$ need not be strict when $\RE_1$ is known, and in fact we would expect the algorithm to perform well with very modest values for $h$ and $\mu$ (see \cite{Gesho_Olson_Titi_2015}).

In addition, the complexity of small viscosity flows requires more computational resources to accurately simulate, but our results indicate that if one has (coarse) measurement data collected on such a flow continuously over a time interval $[0,T]$ it may be possible to construct an accurate computational simulation of the flow over the same time interval, using a much larger value for the viscosity, saving computational resources. Note that one would still need to use the true, smaller viscosity in our simulations after time $T$ to accurately \emph{predict} the behavior of the flow, because we have no data after time $T$. In Section~\ref{sec:subgrid} we test the effectiveness of such an approach numerically.

Lastly, although our primary purpose in the preceding section was to obtain an upper bound on the data assimilation error, in doing so we have obtained a lower bound on the viscosity error, $|\nu_2 - \nu_1|$, in terms of the resulting data assimilation error. In light of this fact, in Section~\ref{sec:param-est} we construct a rudimentary algorithm to estimate the value of the true (but unknown) viscosity, $\nu_1$, using data collected on the flow, $\bu$, and the solution of the data assimilation algorithm, $\bv$. We then test the algorithm numerically.

\subsection{Computational Setting}
All of the following computations were performed on the supercomputer Karst at Indiana University, using \textit{dedalus}, an open source pseudo-spectral python package, available at \url{http://dedalus-project.org}. A $512^2$ computational resolution was used, with a 3/2 dealiasing factor. 
A simple explicit/implicit time stepping scheme was used for each simulation, where the linear terms were handled implicitly, and the nonlinear terms explicitly. 
The spatial domain we consider is $[0,2\pi]$, so $L = 2\pi$. For simplicity, and to limit our assumptions about prior knowledge of the reference solution, we will take the typical velocity to be $1/L$, so that $\RE_1 = \nu_1^{-1}$ and $\RE_2 = \nu_2^{-1}$ in our calculations in the previous sections. With this choice of typical velocity, the Reynolds numbers we define do not characterize the resulting flows in the typical way, so we will instead use the viscosities.

\subsubsection{Reference Solution}
We take our reference solution to be the solution, $\bu^*$, of \eqref{NSEpre} with 
\begin{align*}
	\nu_1 = 0.001,
\end{align*}
and 
\begin{align*}
	f = \sum_{9 < |{\bf k}| < 11 } \widehat{f}_{\bf k} e^{ i{\bf k}\cdot \bx},
\end{align*}
and with the initial condition $\bu^*(0) = 0$. 
Each $\widehat{f}_{ {\bf k} }$ is normally distributed, and scaled so that $|f| = 1$.

\begin{figure}
	\begin{tikzpicture}
		\begin{axis}[%
			xlabel={r},
			ylabel={},
			legend entries={$S(r)$,$S(r)/r^{-5/3}$,$S(r)/r^{-3}$,$S(r)/r^{-5}$},
			tick align=outside,
			tick pos=left,
			legend style={draw=white!80.0!black},
			legend pos={south west},
			x grid style={white!69.01960784313725!black},
			y grid style={white!69.01960784313725!black},
			spectrum style,
			ymode=log,
			xmode=log,
		]%
			\input{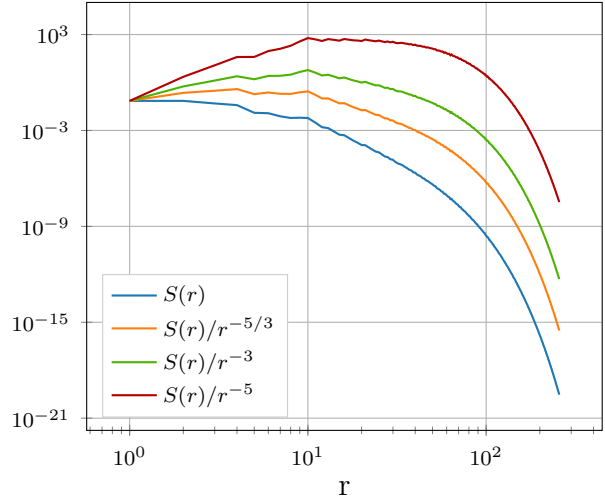}
		\end{axis}
	\end{tikzpicture}
	\caption{Spectrum of the computed reference solution over the time interval \([20,30]\)}\label{fig:spectrum}
\end{figure}

We do not have a closed form solution for $\bu^*$ so instead we approximate it numerically by solving \eqref{NSEpre} computationally over the time interval $[0,30]$. We call the computational approximation we obtain $\bu$, and denote its Fourier transform by $\widehat{\bu}$. So, for all $t\in[0,30]$,
\[ \bu(\bx,t) = \sum_{|{\bf k}|\leq 256} \widehat{\bu}_{\bf k}(t) e^{i{\bf k}\cdot\bx}. \]
Figure~\ref{fig:spectrum} shows the spectrum of $\bu$ over the time interval $[20,30]$, where we define the spectrum, \( S:[0,\infty)\to\nR \), by
\[ S(r) = \frac1{10}\int_{20}^{30} \sum_{ r-\frac12 < |{\bf k}| \leq r+\frac12} |\widehat{\bu}_{\bf k}(t)|^2 dt. \]

\subsubsection{Data Assimilation Parameters}
In the following numerical experiments, we only consider the case that $I_h$ is the projection onto the low modes, i.e. 
\[ I_h(\bu) = (\bx,t)\mapsto\sum_{|{\bf k}|\leq\frac1h} \widehat{\bu}_{\bf k}(t) e^{i{\bf k}\cdot\bx}. \]
We used a spectral method to compute $\bu$, so we can readily construct $I_h(\bu)$. In a practical situation, $I_h(\bu)$ would be given to us and we would have no knowledge of $\bu$; instead, we use $I_h(\bu)$ to compute $\bv$, with the expectation that $\bv(t)\approx\bu(t)$ for all $t$ after a time $t_0$. In Section~\ref{sec:subgrid} we simulate this situation by computing $\bv$ and comparing it to $\bu$.

Before we can compute $\bv$, we will need to choose values for $\mu$ and $h$. The rigorous estimates we have obtained thus far are sufficient conditions, and do not determine the most efficient values of $\mu$ and $h$ in practice. Specifically, for the reference solution we have computed, \( G_1 = 10^6 \), so to satisfy the requirements of Theorem~\ref{thmMainResult1}, we would need \( \mu \sim 10^{12} \) and $h\sim 10^{-6}$. To compute $I_h(\bu)$ with $h = 10^{-6}$, in addition to requiring a large amount of data in practice, would require we increase the computational resolution at least to $200,000^2$. Fortunately, the algorithm works with much less data, and with much smaller $\mu$.

For simplicity, we will only consider \[ \mu = 20, \quad  h = \frac{1}{32}.\]

\subsection{Subgrid Simulations}\label{sec:subgrid}
We are now ready to test the performance of the data assimilation algorithm when $\RE_2\neq\RE_1$. 
We compute the solutions of \eqref{NSEmodified} corresponding to several values of $\nu_2$, with percentage error, $|\nu_2 - \nu_1|/\nu_1$, ranging from $1000\%$ to $0.1\%$.
Each solution is computed over the time interval $[20,30]$ with the initial condition $\bv(20) \equiv 0$. Starting the data assimilation simulation at time $t=20$ is sufficient in this case to ensure that $\bu$ is past a transient (and so is approximating a physical flow), and will be nontrivial at $t=20$ (and therefore differs from $\bv$ at the start of the simulation). 

\begin{figure}
	\begin{tikzpicture}[baseline, trim axis left]
		\begin{axis}[%
			width=.9\textwidth,
			xlabel={$t$},
			ylabel={$|\bu(t) - \bv(t)|$},
			legend entries={$\nu_2 = 0.0008$,$\nu_2 = 0.0009$,$\nu_2 = 0.00099$,$\nu_2 = 0.000999$,$\nu_2 = 0.0009999$,$\nu_2 = 0.0010001$,$\nu_2 = 0.001001$,$\nu_2 = 0.00101$,$\nu_2 = 0.0011$,$\nu_2 = 0.002$,$\nu_2 = 0.011$},
			ymode=log,
			error style,
			]%
			\input{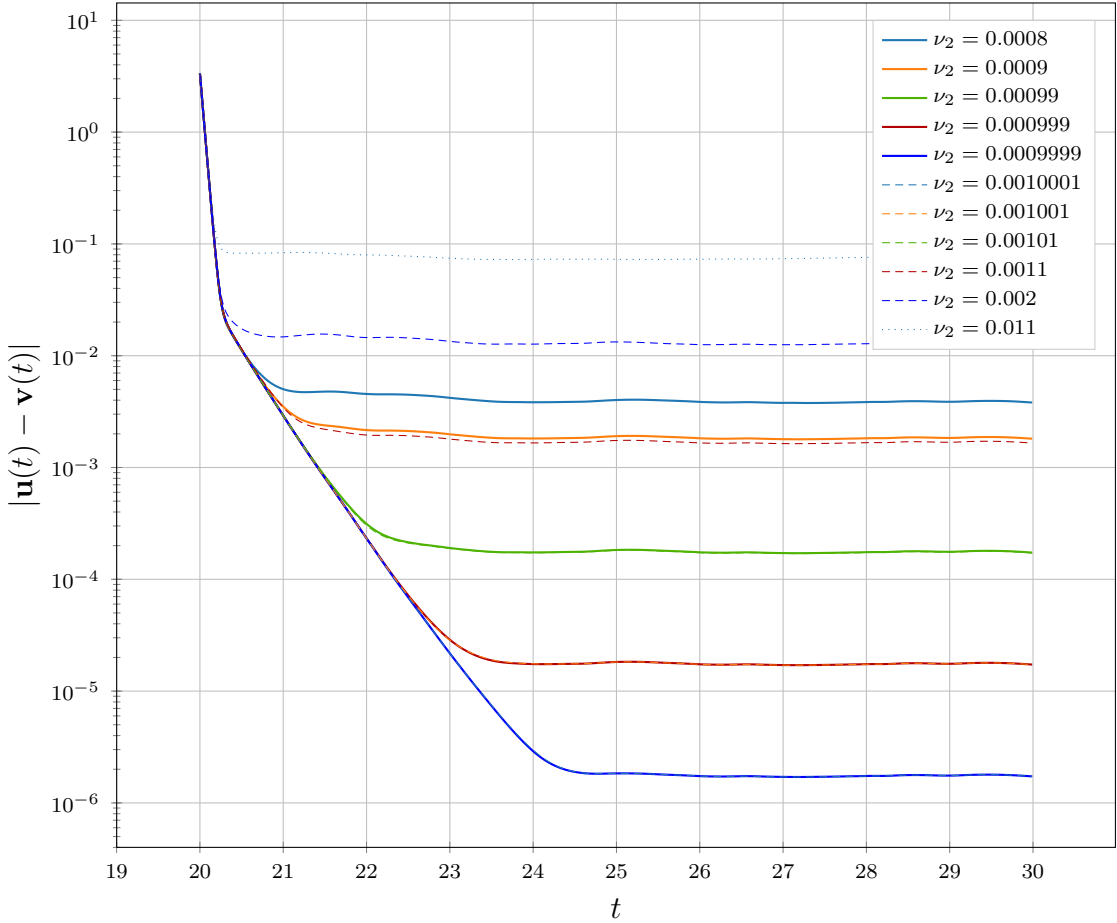}
		\end{axis}
	\end{tikzpicture}
	\caption{The evolution of the $L^2$ error is shown for the solutions of the data assimilation system corresponding to several different values of $\nu_2$. The minimum $L^2$ error achieved decreases as the viscosity error decreases.}\label{fig:error-0}
\end{figure}

Figure~\ref{fig:error-0} shows the resulting $L^2$ error we observe for each simulation when we compare to $\bu$ over the same time interval. We see that for each simulation, after a transient period of fast convergence, the error decreases exponentially at a nearly constant rate before reaching a minimum value. Also, the rates of convergence are the same for each simulation. 

\subsection{Parameter Recovery}\label{sec:param-est}
We can see in Figure~\ref{fig:error-0} that the error in the viscosity value is directly correlated with the minimum error achieved by the corresponding data assimilation solution. This observation motivates the following: given the data $I_h(\bu)$, we can compute $\bv$ and use the minimum error we observe to estimate the true viscosity, $\nu_1$. 

\begin{figure}
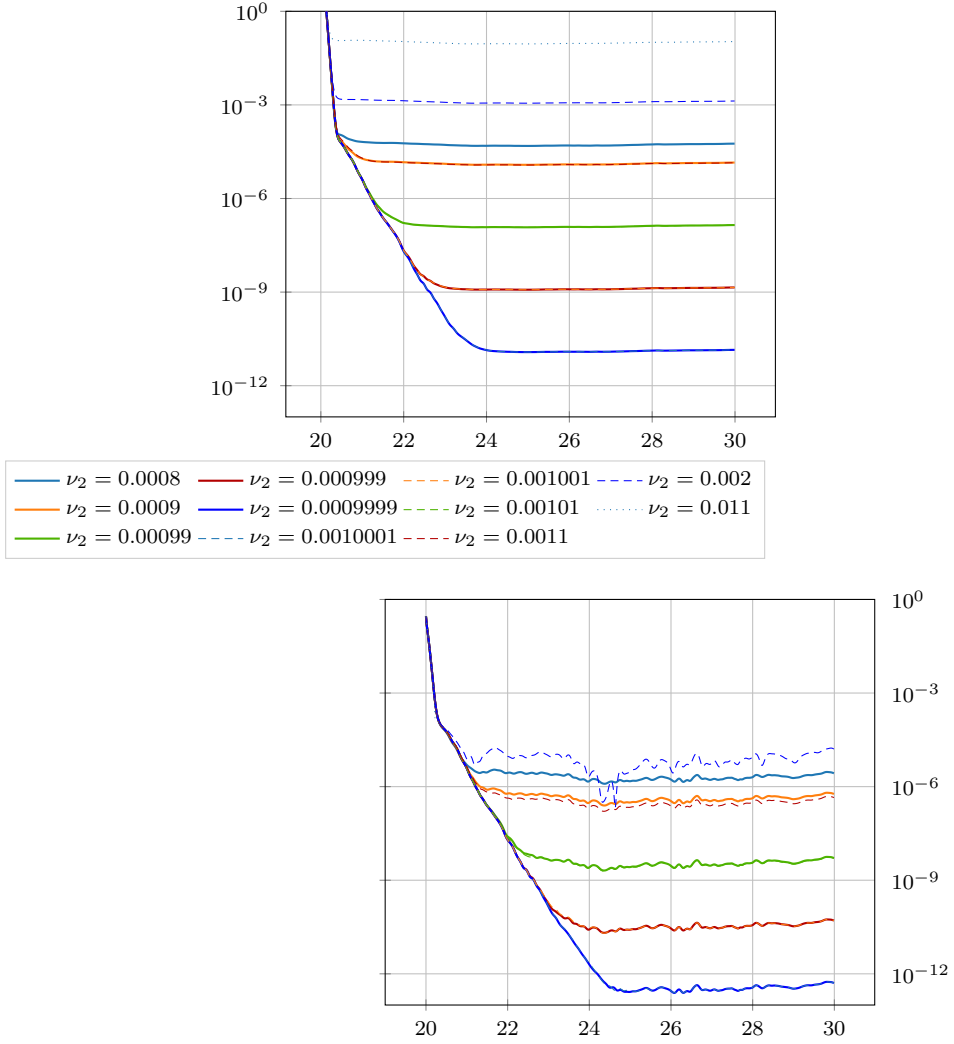

	\begin{tikzpicture}[baseline, trim axis right]
		\begin{axis}[%
			width=.49\textwidth,
			ymode=log,
			ymax={1},
			ymin={1e-13},
			ytickten={3,0,-3,-6,-9,-12},
			error style,
			]%
			\input{Iherror.tex}
		\end{axis}
	\end{tikzpicture}
	\begin{tikzpicture}[baseline, trim axis left]
		\begin{axis}[%
			width=.49\textwidth,
			ymode=log,
			yticklabel pos=upper,
			ytick pos=upper,
			ymax={1},
			ymin={1e-13},
			ytickten={3,0,-3,-6,-9,-12},
			error style,
			legend entries={$\nu_2 = 0.0008$,$\nu_2 = 0.0009$,$\nu_2 = 0.00099$,$\nu_2 = 0.000999$,$\nu_2 = 0.0009999$,$\nu_2 = 0.0010001$,$\nu_2 = 0.001001$,$\nu_2 = 0.00101$,$\nu_2 = 0.0011$,$\nu_2 = 0.002$,$\nu_2 = 0.011$},
			legend style={at={(0,1.1)},anchor=south},
			transpose legend,
			legend columns=3,
			]%
			\input{RHS.tex}
		\end{axis}
	\end{tikzpicture}
	\caption{Shown on the left is $\mu|I_h(\bv) - I_h(\bu)|^2$ vs time for several different values of $\nu_2$. We see that in each case, $|I_h(\bv) - I_h(\bu)|$ reaches a minimum value, which is smaller for $\nu_2$ closer to $\nu_1$. On the right is the value of the right hand side of \eqref{NSE-modiff-v-Ih} for the same values of $\nu_2$. We see that the values on the right are negligible compared to the error values on the left. }\label{fig:Ih-error}
\end{figure}

Although $I_h(\bu)$ is sufficient to compute $\bv$, we would need to have $\bu$ to compute $|\bu - \bv|$. Fortunately, we see that $|\nu_2 - \nu_1|$ and $|I_h(\bu) - I_h(\bv)|$ are also correlated, as can be seen in Figure~\ref{fig:Ih-error}.
With this in mind, we will now study this correlation, so that, once it's nature is established, we can use it to develop an algorithm to estimate $\nu_1$.

\subsubsection{A posteriori Error Estimate}\label{secAposteriori}
The result in Theorem~\ref{thmMainResult1}, in addition to being in terms of the true error (as opposed to the error of only the interpolations of $\bu$ and $\bv$), establishes bounds for the data assimilation error in terms of the Grashof number. We are now considering a situation where we have access to $\bv$, and so would like to obtain a sharper estimate on the error by allowing it to be in terms of $\bv$ instead of $G$.

As in the proof of Theorem~\ref{thmMainResult1}, let $\bw = \bv-\bu$. Using the facts that
\[ B(\bw,\bu)+ B(\bv, \bw) = B(\bu,\bw) + B(\bw,\bv),\]
and 
\[ (\RE_2^{-1} - \RE_1^{-1}) A\bu + \RE_2^{-1} A\bw = (\RE_2^{-1} - \RE_1^{-1}) A\bv + \RE_1^{-1} A\bw, \]
we can replace \eqref{NSE_modiff} with
\begin{align*}
	\bw_t + B(\bw,\bv)+ B(\bu, \bw) = (\nu_2-\nu_1) A \bv + \nu_1 A \bw - \mu P_\sigma( I_h(\bw)).
\end{align*}

Now, we apply $I_h$ to both sides of this equation and obtain
\begin{align*}
	\partial_t I_h(\bw) + I_h(B(\bw,\bv)+ B(\bu, \bw)) \\
	= (\nu_2-\nu_1) I_h(A \bv) + \nu_1 I_h(A \bw) - \mu I_h(P_\sigma( I_h(\bw))).
\end{align*}
Next, we take the inner product with \(I_h(\bw)\) and use the fact that
\begin{equation*}
	- \mu\ip{I_h(P_\sigma( I_h(\bw) ))}{I_h(\bw)} = - \mu |I_h(\bw)|^2. 
\end{equation*}
The result, after rearranging terms, is
\begin{align}\label{eqn:param-diff}
	\frac12\frac{d}{dt} |I_h(\bw)|^2 + (\nu_1-\nu_2) \ip{I_h(A \bv)}{I_h(\bw)} + \mu |I_h(\bw)|^2 
	\\
	= \ip{I_h\left( \nu_1 A \bw - B(\bw,\bv) - B(\bu, \bw)\right)}{I_h(\bw)}. 
\end{align}

We have observed in each of our simulations that there is a time at which the error $|I_h(\bw)|$ reaches a minimum value and thereafter remains constant; then $\frac{d}{dt} |I_h(\bw)| \approx 0$, and the above equation reduces to
\begin{align}\label{NSE-modiff-v-Ih}
	(\nu_1-\nu_2) \ip{I_h(A \bv)}{I_h(\bw)} + \mu |I_h(\bw)|^2 
	\\\notag
	= \ip{I_h\left( \nu_1 A \bw - B(\bw,\bv) - B(\bu, \bw)\right)}{I_h(\bw)}. 
\end{align}

Note that all of the terms on the left hand side of \eqref{NSE-modiff-v-Ih} except $\nu_1$ are explicitly computable from data observations. However, on the right hand side, one would need $\bu$ to compute $B(\bw,\bv)$ and $B(\bu,\bw)$. Also, although in the periodic setting, $A$ commutes with the projection onto the low Fourier modes, $A$ might not commute with other types of interpolation operators $I_h$, in which case one could not compute $I_h(A\bw)$ exactly from the observations $I_h(\bu)$. 

However, we note that in terms of units, each of the terms in \eqref{NSE-modiff-v-Ih} decreases quadratically with $\bw$ as $\bw\to0$ (with the exception of $(\nu_1-\nu_2) \ip{I_h(A \bv)}{I_h(\bw)}$), but we control $\mu$ and have chosen $\mu$ large enough that $\mu|I_h(\bw)|^2$ dominates the terms on the right hand side, as can be seen in Figure~\ref{fig:Ih-error}. Therefore, we propose an approximation formed by dropping these terms from the equation, and solving (approximately) for $\nu_1$, thereby obtaining
\begin{equation}\label{param-est-main}
\nu_1 \approx \nu_2 - \mu \frac{|I_h(\bw)|^2}{ \ip{I_h(A \bv)}{I_h(\bw)} }.
\end{equation}
Since each time on the right-hand side now depends only on given or observable quantitesm, 
This approximation motivates an iterative scheme for recovering the viscosity.  We therefore test \eqref{param-est-main} as a means of recovering $\nu_1$, using the data from our simulations. We obtain the approximation $\tilde{\nu}_1$ iteratively, using \eqref{param-est-main} for each of the simulations performed in Section~\ref{sec:subgrid} at time $t=24$, and compare to $\nu_1$. The results are shown in Table~\ref{tbl:1}. In each case, \eqref{param-est-main} produces a much better approximation of the true $\nu_1$, showing at least an $80\%$ improvement. 

\begin{table}[ht]
	\centering
	\caption{ }\label{tbl:1} 
	\begin{tabular}{|c||c|c|c|c|c|}\hline
		\( \nu_2 \) &  
		\( |I_h(\bw)|^2 \) & 
		\( \ip{I_h(A \bv)}{I_h(\bw)} \) &
		\( \tilde{\nu}_1 \) &
		\( |\tilde{\nu}_1 - \nu_1| \) &
		\( \frac{|\tilde{\nu}_1 - \nu_1|}{|\nu_2-\nu_1|} \) 
		\\ \hline \hline
				0.00080000 & 1.564e-03 & -2.462e-01 & 9.986e-04 & 1.381e-06 & 0.69\% \\ \hline
		0.00090000 & 7.786e-04 & -1.226e-01 & 9.988e-04 & 1.155e-06 & 1.15\% \\ \hline
		0.00099000 & 7.760e-05 & -1.223e-02 & 9.999e-04 & 1.495e-07 & 1.49\% \\ \hline
		0.00099900 & 7.762e-06 & -1.222e-03 & 1.000e-03 & 1.423e-08 & 1.42\% \\ \hline
		0.00099990 & 8.309e-07 & -1.223e-04 & 1.000e-03 & 1.288e-08 & 12.88\% \\ \hline
		0.00100010 & 8.339e-07 & 1.221e-04 & 1.000e-03 & 1.391e-08 & 13.91\% \\ \hline
		0.00100100 & 7.765e-06 & 1.222e-03 & 1.000e-03 & 1.328e-08 & 1.33\% \\ \hline
		0.00101000 & 7.755e-05 & 1.222e-02 & 1.000e-03 & 1.551e-07 & 1.55\% \\ \hline
		0.00110000 & 7.732e-04 & 1.218e-01 & 1.002e-03 & 1.826e-06 & 1.83\% \\ \hline
		0.00200000 & 7.570e-03 & 1.183e+00 & 1.031e-03 & 3.121e-05 & 3.12\% \\ \hline
		0.01100000 & 6.750e-02 & 9.430e+00 & 1.336e-03 & 3.356e-04 & 3.36\% \\ \hline

	\end{tabular}
\end{table}

\newcommand{\IP}[2]{\int_s^t \ip{I_h(#1(\tau))}{I_h(#2(\tau))}d\tau}
\newcommand{\NORM}[1]{\int_s^t |I_h(#1(\tau))|^2 d\tau}
To avoid waiting until the time derivative of the error becomes negligible, or to include the possibility that the time derivative has non-negligible oscillations, we can choose to leave the time derivative in \eqref{eqn:param-diff}. Let $t > s \geq t_0$. We then integrate \eqref{eqn:param-diff} over the time interval $[s,t]$, to obtain
\begin{align*}\label{NSE-modiff-v-Ih-2}
	\frac12 |I_h(\bw(t))|^2 - \frac12 |I_h(\bw(s))|^2 + (\nu_1-\nu_2) \IP{A \bv}{\bw} + \mu \NORM{\bw} 
	\\
	= \IP{\left( \nu_1 A \bw - B(\bw,\bv) - B(\bu, \bw)\right)}{\bw}. 
\end{align*}
Then, dropping the terms on the right hand side as before and solving for $\nu_1$, we obtain
\begin{equation}\label{param-est-main-2}
	\nu_1 \approx \nu_2 - \frac{\mu\NORM{\bw} + \frac12 |I_h(\bw(t))|^2 - \frac12|I_h(\bw(s))|^2  }{ \IP{A \bv}{\bw} }.
\end{equation}

\subsubsection{Algorithms}
We next use \eqref{param-est-main} and \eqref{param-est-main-2} to devise algorithms capable of recovering $\nu_1$ using only the data $I_h(\bu)$ over a time interval $[t_0,T]$. The first algorithm (Algorithm~\ref{alg:1}) utilizes \eqref{param-est-main}. Algorithm~\ref{alg:2} describes a method to recover $\nu_1$ using \eqref{param-est-main-2} instead of \eqref{param-est-main}.
\begin{algorithm}[ht]
\caption{}\label{alg:1}
\def\Compute{\State \textbf{compute }}
\def\Choose{\State \textbf{choose }}
\def\Input{\State \textbf{input }}
\begin{algorithmic}[0]
	\small
	\Input $I_h(\bu)$ on $[t_0,T]$ \Comment{available reference solution data}
	\Input $\nu_2$ \Comment{an initial estimate for $\nu_1$}
	\Input $dt > 0$ \Comment{time step}
	\Input $\epsilon > 0$ \Comment{tolerance for machine precision}
	\Input$\delta \in (0,1)$ \Comment{tolerance for convergence}

	\State $t \gets t_0$
	\State $\bv(t_0) \gets 0$

	\While{$|I_h(\bu(t)) - I_h(\bv(t))| > \epsilon$ \and $t<T$}
		\Compute $\bv(t+dt)$ using viscosity $\nu_2$ and feedback $I_h(\bu(t))$
		\If{$|I_h(\bu(t+dt)) - I_h(\bv(t+dt))| \geq (1-\delta)|I_h(\bu(t))-I_h(\bv(t))|$} 
			\If{$|I_h(\bu(t+dt)) - I_h(\bv(t+dt))| < |I_h(\bu(t_0))-I_h(\bv(t_0))|$} 
				\Compute $\tilde{\nu}_1$ using \eqref{param-est-main} at time $t+dt$
				\State $t_0 \gets t+dt$
				\State $\nu_2 \gets \tilde{\nu}_1$
			\Else
				\State \Return $\nu_2$
			\EndIf
		\EndIf
		\State $t\gets t+dt$
	\EndWhile

	\State \Return $\nu_2$
\end{algorithmic}
\end{algorithm}

\begin{algorithm}[ht]
\caption{}\label{alg:2}
\def\Compute{\State \textbf{compute }}
\def\Choose{\State \textbf{choose }}
\def\Input{\State \textbf{input }}
\begin{algorithmic}[0]
	\small
	\Input $I_h(\bu)$ on $[t_0,T]$ \Comment{available reference solution data}
	\Input $\nu_2$ \Comment{an initial estimate for $\nu_1$}
	\Input $dt > 0$ \Comment{time step}
	\Input $I > 0$ \Comment{time to wait before computing time averages}
	\Input $J > 0$ \Comment{length of time interval used to compute time averages}
	\Input $\epsilon > 0$ \Comment{tolerance for machine precision}

	\State $\bv_0 \gets 0$
	\State $t \gets t_0$

	\While{$|I_h(\bu(t)) - I_h(\bv(t))| > \epsilon$ \and $t_0 + I + J < T$}
		\Compute $\bv(t)$ on $[t_0,t_0+I+J]$ using viscosity $\nu_2$, feedback $I_h(\bu(t))$, IC $\bv(t_0) = \bv_0$, and time step $dt$.
		\Compute $\tilde{\nu}_1$ using \eqref{param-est-main-2} over the time interval $[t_0+I, t_0+I+J]$.
		\State $\bv_0 \gets \bv(t_0 + I + J)$
		\State $t_0 \gets t_0 + I + J$

		\State $\nu_2 \gets \tilde{\nu}_1$
	\EndWhile

	\State \Return $\nu_2$
\end{algorithmic}
\end{algorithm}

\begin{figure}
	\begin{tikzpicture}[baseline, trim axis left]
		\begin{axis}[%
			width=.9\textwidth,
			xlabel={$t$},
			legend entries={Alg~1 $|\bu(t) - \bv(t)|$,Alg~1 $|I_h(\bu(t)) - I_h(\bv(t))|$,Alg~1 $|\nu_2 - \nu_1|/\nu_1$,Alg~2 $|\bu(t) - \bv(t)|$,Alg~2 $|I_h(\bu(t)) - I_h(\bv(t))|$,Alg~2 $|\nu_2 - \nu_1|/\nu_1$},
			ymode=log,
			PE style,
			]%
			\input{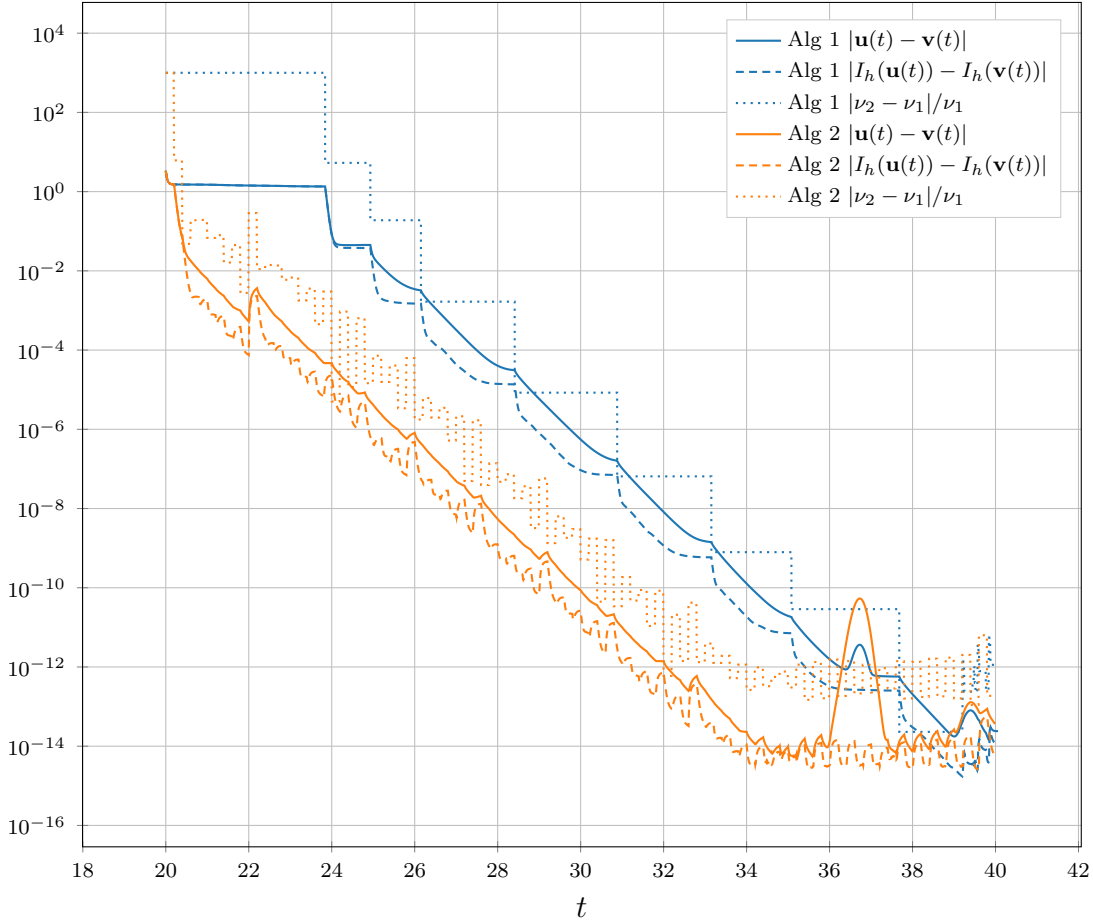}
		\end{axis}
	\end{tikzpicture}
	\caption{The evolution of the $L^2$ error is shown for the solutions of the data assimilation systems corresponding to Algorithm~\ref{alg:1} and Algorithm~\ref{alg:2}, as well as that of the relative error in the approximate viscosity. The $\nu_k$ are chosen and the equations updated following the procedures outlined in the Algorithms.}\label{fig:PE-alg-2}
\end{figure}

Figure~\ref{fig:PE-alg-2} shows the errors observed during the process of applying Algorithm~\ref{alg:1} and Algorithm~\ref{alg:2} to our reference solution.
Algorithm~\ref{alg:1}, tested with our reference solution and an initial guess of $\nu_2 = 1$, after $10$ iterations produced an end approximation of $\tilde{\nu}_1 = 0.00100000000000113301$ (an absolute error of $\approx1.133\times10^{-15}$). With similar performance, Algorithm~\ref{alg:2}, tested under the same conditions, after $100$ iterations produced an end approximation of $\tilde{\nu}_1 = 0.00099999999999981332$ (an absolute error of $\approx1.867\times10^{-16}$). In both cases, the results are accurate to within machine precision.

\section{Conclusion}
In this article, we presented and analyzed a new way to recover unknown parameters of a system (in this case, the Reynolds number, or equivalently, the viscosity), using a continuous data assimilation approach for the 2D incompressible Navier-Stokes equations.  This means that even in the case where the viscosity is unknown and one only has sparse observational data, one may still obtain convergence to the true solution by using the AOT algorithm in combinatoin with the algorithms proposed here.  In addition, our new algorithms allow one to update the viscosity in real time using only observational data, and we showed computationally that the true solution and the true viscosity are recovered to within machine precision, exponentially fast in time.  An analytical proof of this will be the subject of a forthcoming work, which will also explore the extension of the algorithm to other physical systems.

In addition, and as a lead-up to the new algorithms, we proved analytically that in the case of an inaccurately known viscosity, the large-time error produced by the AOT algorithm is controlled by the error in the viscosity.  

Since our new algorithms involve changing the viscosity mid-simulation, we also examined the corresponding viscosity sensitivity equations.  Specifically, we proved the existence and uniqueness of global solutions to these equations.  A byproduct of the proof is that the sensitivity of solutions to the equations involved in the algorithm are bounded in appropriate spaces.  Hence, changing the viscosity mid-simulation does not result in major aberrations in the solution.  We note that in the present context, our proof is somewhat non-standard, in that we proved the existence by showing that the difference quotients converge (or at least, have a subsequence that converges) to a solution of the equations.  We also note that this appears to be the first such rigorous proof that the sensitivity equations for the 2D Navier-Stokes equations have a unique solution, although formal proofs have been given in other works, cited above.

%
%
%

\section*{Acknowledgements}
 \noindent
 The research of E.C. was supported in part by the NSF GRFP grant no. 1610400. The research of A.L. was supported in part by the NSF grant no. DMS-1716801. The research of J.H. was supported in part by NSF Grant DMS-1517027. Computational resources were provided by Lilly Endowment, Inc., through its support for the Indiana University Pervasive Technology Institute, and in part by the Indiana METACyt Initiative. The Indiana METACyt Initiative at IU was also supported in part by Lilly Endowment, Inc.
 
\bibliographystyle{abbrv}

\end{document}